\documentstyle[12pt,fleqn,leqno]{article}
\newtheorem{thm}{Theorem}[section]
\newcommand{\be}{\begin{equation}}
\newcommand{\ee}{\end{equation}}
\newcommand{\ba}{\begin{array}}
\newcommand{\ea}{\end{array}}

\renewcommand{\a}{\alpha}
\newcommand{\kp}{\kappa}
\renewcommand{\b}{\beta}
\renewcommand{\l}{\lambda}
\renewcommand{\t}{\theta}

\newcommand{\mg}{\rm}
\renewcommand{\em}{\it}
\newcommand{\bea}{\begin{eqnarray}}
\newcommand{\eea}{\end{eqnarray}}

\newcommand{\G}{\Gamma}
\newcommand{\g}{\gamma}

\begin{document}
\newtheorem{lem}[thm]{Lemma}
\newtheorem{cor}[thm]{Corollary}
\title{\bf Generalized Orthogonality and Continued
Fractions\thanks{
Research partially supported  by NSF grant DMS 9203659 and NSERC
grant
 A5384}}
\author{Mourad E. H.  Ismail and David R. Masson }
\date{July 1994}
\maketitle
\begin{abstract}
 The connection between continued fractions and orthogonality which
 is familiar for $J$-fractions and $T$-fractions is extended to 
 what we call $R$-fractions of type I and II. These continued
fractions
 are associated with recurrence relations that correspond to
 multipoint rational interpolants. A Favard type theorem is proved 
 for each type. 
 We then study explicit models which lead to
biorthogonal rational functions.
\end{abstract}

\bigskip 
{\bf Running title}:$\;$  Continued R-fractions

\bigskip
{\em 1990  Mathematics Subject Classification}:  Primary 40A15,
41A20
Secondary 33C45, 33D45.

{\em Key words and phrases}. 
Continued fractions,  $T$-fractions, $R$-fractions,
integral representations,
Pincherle's theorem, biorthogonal polynomials and rational
functions,
contiguous relations, Cauchy beta integral, $q$-beta integrals,
 Askey-Wilson polynomials.

\vfill\eject
\setcounter{section}{1}

{\bf 1. Introduction}. 
There is a close connection between the theory of
orthogonal polynomials and continued fractions. This connection is 
clearest in the case of the classical moment problem where
orthogonality 
is with respect to a positive measure whose Stieltjes transform
is represented by a positive definite J-fraction (see
\cite{Jo:Th1},
\cite{Wa}). This connection is known to extend to the
quasi-definite 
case where the measure and the J-fraction are no longer 
positive and positive definite, respectively \cite{Ch},
\cite{Gu:Is}.
Another
extension is by means of T-fractions and their connection with the 
trigonometric moment problem \cite{Jo:Nj:Th}, \cite{Jo:Th1},
\cite{Sz}
. With
both types of extensions one has a Favard type theorem
establishing orthogonality with respect to a unique
moment functional, given the three term recurrence 
relation satisfied by the polynomials under consideration. More
importantly there is always a spectral measure(s) associated
with these problems and by determining the regions where the 
associated continued fraction diverges one can locate the support
of the corresponding spectral measure(s). Furthermore by finding
the limits to which the associated continued fraction converges 
one can compute, at least in theory, the spectral measure(s).

In this paper we shall examine two further natural extensions
which, from a different point of view, are associated with
multipoint
rational interpolants or Pad\'e approximants (\cite{Go},
\cite{Go:Lo},
\cite{Nj}, \cite{He:Nj},
\cite{St:To}).
We call the new types of continued fractions "$R$-fractions".
We establish a Favard type theorem for $R$-fractions. We also
establish the existence of a natural Borel measure associated with
$R$-fractions. In general we show the existence of 
  rational biorthogonality 
 with respect to this 
 Borel measure. We demonstrate by
 means of examples how our 
general approach can yield explicit biorthogonal systems of
rational functions and can lead to the evaluation
of some new types of beta integrals.
 
 In Section 2 we consider what we call continued fractions of 
 type $R_I$. Such continued fractions are associated with the 
 polynomial recurrence relation
 \begin{eqnarray}
 P_n(x) - (x - c_n)P_{n-1}(x) + \lambda_n (x- a_n) P_{n-2}(x)
 = 0.
 \end{eqnarray}
 In addition to proving a Favard type theorem for $R_I$-fractions
in
 Section 2, we also establish a biorthogonality relation for 
  rational functions. Section 2 contains, as an application of
  our results, a new evaluation of a $q$-beta integral of Ramanujan
and
  an alternate derivation of a biorthogonality relation due to
Pastro 
  \cite{Pa}. Pastro's original proof involves the aforementioned
integral
  of Ramanujan. It is worth noting that we obtain both the
evaluation
   of the $q$-beta integral and the biorthogonality relation from
the
   same result on $R_I$-fractions. An additional result is the
evaluation 
   of the Herglotz transform of the Borel measure involved.

 In Section 3 we consider continued fractions of type $R_{II}$ 
 which are associated with polynomial recurrence relation 
 \begin{eqnarray}
 P_n(x) - (x - c_n)P_{n-1}(x) + \lambda_n \, (x - a_n)(x - b_n)
 P_{n-2}(x) = 0.
 \end{eqnarray}
 In Section 3 we follow the same plan  as in Section 2. We prove a 
 representation theorem for moment functionals associated with
 (1.2) and study the convergence of the corresponding continued
 fraction. We also identify biorthogonal rational functions that
 arise naturally from the recurrence relation (1.2). 

 The letter $R$ in our terminology is used to emphasize that the  
  corresponding continued fractions are associated with rational
  interpolation
 and rational biorthogonality.
 
 To illustrate our ideas we include seven examples where we have
 applied our general results to concrete situations. We have
already
 mentioned that in Section 2, we rederive a result of Pastro
\cite{Pa}. In
 Section 3 we give two new systems of biorthogonal rational
functions
  associated with continued fraction of type $R_{II}$. 
The first is based on the $R_{II}$  recurrence relation with
constant coefficients which corresponds to  Chebyshev polynomials.
The second 
uses $_2F_1$ contiguous relations and yields a weight function
which is the 
integrand in the Cauchy beta integral. Section 4
  contains two  examples of biorthogonal ${}_4\phi_3$  
  rational functions based on solutions to the Askey-Wilson three
term 
recurrence  relation.  The first rederives a system of biorthogonal
rational
 functions due to Al-Salam and Ismail \cite {Al:Is}. The second
turns out to 
yield special cases of  a biorthogonality related to
$q^{-1}$-Hermite
 polynomials, 
\cite {Is:Ma}. Section 5 contains two additional examples based on
a modification of the Askey-Wilson recurrence. The first is new and
unexpected.  
It provides  an integral representation for a  ${}_2\phi_1$ basic
 hypergeometric function and  a special case of it leads to
rational
 biorthogonality on $[-1, 1]$ based on the elementary integral
\be
\frac{2}{\pi} \int_{-1}^1 \frac{\sqrt{1-x^2}\, dx}{(1 - 2 \a x +
\a^2)(1 - 2\delta x +\delta^2)} = \frac{1}{1 - \a\delta}.
\ee
The second, which is originally due to M. Rahman \cite {Ra},
evaluates a new $q$-beta integral of the type studied by Askey
 and Wilson and derives an associated system of  rational functions
biorthogonal
on $[-1, 1]$. The biorthogonality is between a
 $\{{}_4\phi_3\}$ and  $\{{}_6\phi_5\}$ pair.
The methodology followed in the examples combines contiguous
relations 
for hypergeometric or basic hypergeometric functions with 
Pincherle's theorem \cite{Jo:Th1}, \cite{Lo:Wa}.
This methodology proved fruitful
in other types of continued fractions, \cite{Is:Li}, \cite{Gu:Is},
 \cite{Is:Le:Va:Wi}, \cite{Ma}, \cite{Wi1}, \cite{Wi2}.

We shall mostly follow the notation of Gasper and Rahman in 
\cite{Ga:Ra}. The shifted  factorials and multishifted 
factorials are
\be
(a)_0 := 1, \quad (a)_n:= \prod_{j=1}^{n}(a + j -1), \quad n > 0,
\ee
and
\be
(a_1, a_2, \cdots, a_k)_n := \prod_{j=1}^k (a_j)_n,
\ee
respectively, 
and a hypergeometric function is
\bea
{}_rF_s(a_1,a_2, \cdots, a_r; b_1, b_2, \cdots, b_s; z) =
{}_rF_s\left(
\left. \ba{c} a_1, a_2, \cdots, a_r \\ b_1, b_2, \cdots, b_s
\ea \right| z \right)
\eea
 \bea
 \quad \quad  := \sum_{n=0}^{\infty} \frac{(a_1, a_2, \cdots,
a_r)_n}
 {(b_1, b_2,\cdots, b_s)_n}\, \frac{z^n}{n!}. \nonumber
\eea
The $q$-shifted and multishifted factorials are similarly defined
by
\be
(a;q)_0:=1,\quad (a;q)_n := \prod_{j = 1}^n (1-a q^{j-1}),\quad  n
> 0,
\; or \; n = \infty,
\ee
 \be
 (a_1, a_2, \cdots, a_k; q)_n := \prod_{j = 1}^k (a_j;q)_n.
 \ee
 A basic hypergeometric series is
 \bea
 {}_r\phi_s(a_1, a_2, \cdots, a_r; b_1, b_2, \cdots; q, z) 
 = {}_r\phi_s \left( \left. \ba{c}a_1, a_2, \cdots, a_r \\
 b_1, b_2, \cdots, b_s \ea \right| q, z\right)
 \eea
 \bea
 \quad \quad :=\sum_{n=0}^\infty \frac{(a_1, a_2, \cdots, a_r;
q)_n}
 {(b_1, b_2, \cdots, b_s; q)_n} \; \frac{[(-1)^n q^{n(n-1)/2}]^{1 
 + s -r} z^n}{(q; q)_n}. \nonumber
\eea
 Another useful notation is the short hand notation for a very 
 well-poised basic series. We shall use
 \bea
 {}_8W_7(a; b, c, d, e, f; z) :={}_8\phi_7\left( \left. \ba{c}
 a, \sqrt{aq},\; -\sqrt{aq},\; b, \;c,  \; d, \;e, \;f \\ \sqrt{a},
-\sqrt{a},
 aq/b, aq/c, aq/d, aq/e, aq/f \ea \right|q, z\right).
 \eea

The biorthogonality based on (1.3) which we alluded to earlier  is
\be
\frac{2}{\pi} \int_{-1}^1 f_m(x; \delta,  \a)f_n(x; \a, \delta)
\sqrt{1-x^2}\, dx = \frac{(\a\delta)^n\, (q, q; q)_n}{(\a\delta,
\a\delta ;q)_n
(1 - \a\delta q^{2n})} \;\delta_{m,n},
\ee
with $max\{|\a|,\; |\delta|\} < 1$ and $f_n$ is defined by
\be
f_n(x; \a, \delta) := \frac{1}{1- 2\a x + \a^2} {}_4\phi_3\left(
\left. \ba{c}
 q^{-n}, \a\delta q^{n}, \a e^{i\t}, \a e^{-i\t}\\
 \a\delta,\;  q\a e^{i\t}, \; q\a e^{-i\t} \ea \right| q, q\right),
 \; x = \cos \t.
\ee
The functions in (1.12) can be expressed in the simpler form
\be
f_n(x; \a, \delta) = \sum_{k = 0}^n \frac{(q; q)_n}{(q; q)_k(q;
q)_{n-k}}
\frac{(\a\delta q^n; q)_k}{(\a\delta;q)_k}
\frac{(-1)^kq^{k(1-k)/2}}{1 - 2\a q^kx + \a^2q^{2k}}.
\ee
The expression $\frac{(q; q)_n}{(q; q)_k(q; q)_{n-k}}$ is the
Gaussian 
binomial coefficient. The functions $\{f(x; \a, \delta)\}$ are
orthogonal with respect to 
the weight function of the Chebyshev polynomials of the second kind
and may play an important role in a future theory of biorthogonal
rational functions.

\bigskip

\setcounter{section}{2}
 \setcounter{equation}{0}
 
{\bf 2.$R_I$-fractions}. We begin this section with a Favard type
theorem
 for the polynomial recurrence (2.1) below. We then outline how
 orthogonality may be realized in terms of the properties 
of what we will call an $R_I$ type continued fraction.
We also include an example to illustrate this approach. 

We first establish a Favard type theorem which is Theorem 2.1. 
Consider a system of monic orthogonal polynomials generated by
\be
P_n(x) = (x-c_n)P_{n-1}(x) - \lambda_n\,(x-a_n)\, P_{n-2}(x),
\ee
\begin{eqnarray}
P_{-1}(x):= 0, \; P_0(x) := 1, \nonumber
\end{eqnarray}
where
\begin{eqnarray}
\lambda_{n+1} \ne 0, \quad P_n(a_{n+1}) \ne 0, n = 1, 2,
\dots. \nonumber
\end{eqnarray}
The recurrence relation in (2.1) can be renormalized to 
yield the rational recurrence relation
\be
(x-a_{n+1})\, R_n(x) - (x-c_n) R_{n-1}(x) +\, \lambda_n R_{n-2}(x)
 = 0
 \ee
 with the same initial conditions, namely
 \be
 R_0(x) := 1, R_{-1}(x) := 0.
 \ee
 The renormalization is given explicitly by
 \be
 R_n(x) := P_n(x)/\prod_{k = 1}^n(x - a_{k+1}).
 \ee
 \begin{thm}
 Associated with the recurrence relation (2.1) there is a linear 
 functional ${\cal L}$ defined on the span
 of $\{x^kR_n(x)\}_{n,k=0}^\infty$
 mapping it into ${\cal C}$, and normalized by ${\cal L}[1] = 
 \lambda_1 \ne 0$, such that the orthogonality relation
 \be
 {\cal L}[x^kR_n(x)] = 0, 0 \le k < n, 
 \ee
 holds.
 Furthermore the functional values ${\cal L}[x^n]$ and ${\cal L}
 [\prod_{k=1}^n (x - a_{k+1})^{-1}]$,
 $n = 1, 2, \dots$ are uniquely determined in terms of the
sequences
 $\{a_{n+1}, c_n, \l_n\, :\, n = 1, 2, \cdots\}$.
 \end{thm}
 {\bf Proof}. Define a linear functional ${\cal L}$ 
  whose
 action on $R_n(x)$ and $x^{n-1}R_n(x)$ is given by
 \be
{\cal L}[1] :=\lambda_1, \quad\; and \quad\;  {\cal L}[R_n(x)]
= {\cal L}[x^{n-1}R_n(x)]
= 0, for\; n \ge 1.
\ee
We first establish (2.5) for $1 \le k < n-1$. Clearly (2.6)
yields 
\be
0 ={\cal L}[R_1(x)] = {\cal L}[xR_2(x)] = {\cal L}[x^2R_3(x)].
\ee
Apply ${\cal L}$ to (2.2) with $n = 3$ to obtain ${\cal L}[xR_3(x)]
= 0$. Thus (2.5) holds for $1 \le n \le 3$. When $n  > 3$, 
formula (2.5) follows from (2.2) and (2.6) by induction over $n$.
To prove the latter part of theorem observe that ${\cal L}[P_1
(x)/(x-a_2)] = 0$ implies ${\cal L}[1 + P_1(a_2)/(x - a_2)] = 0$.
Using ${\cal L}[1] = \lambda_1$ we then obtain  
\bea
{\cal L} [(x - a_2)^{-1}]  = -\lambda_1/P_1(a_2).  \nonumber
\eea
Similarly from ${\cal L}[P_2(x)/\{(x-a_2)(x-a_3)\}] = 0$ we find
\bea
{\cal L}[1/\{(x-a_2)(x - a_3)\}] = \lambda_1(\lambda_2 + a_3 - a_2)
/\{P_2(a_3)P_1(a_2)\}.   \nonumber
\eea
We next use ${\cal L}[x^k  R_n(x)]$ for $n = 1,2$ and 
$k < n$ to obtain
\bea
{\cal L}[P_2(x)/(x - a_2)] = {\cal L}[-\lambda_2 + (x - c_2)P_1
(x)/(x-a_2)] = {\cal L}[x - c_1 - \lambda_2] = 0. \nonumber
\eea
Thus ${\cal L}[x] = (c_1+\lambda_2)\,\l_1$. We now 
continue by induction on
$n$.  Thus for each new $n$ we use ${\cal L}[R_n(x)] = 0$ and
we obtain a value for ${\cal L}[1/\prod_{k=1}^n (1 - a_{k+1})]$
while from ${\cal L}[x^{n-1}R_n(x)] = 0$ and the recurrence
relation
(2.2) we evaluate ${\cal L}[x^{n-1}]$. This  establishes our
theorem.

 \begin{cor}We have
 \be
 {\cal L}[P_n(x)\,R_n(x)]= {\cal L}[x^n R_n(x)] = 
 \l_1\, \l_2\, \dots \,\l_{n+1}, \; n\ge 0. 
 \ee
\end{cor}
{\bf Proof}. Multiply (2.2) by $x^{n-2}$, apply ${\cal L}$ and
then use
the orthogonality relation (2.5). This yields the two term
recurrence relation ${\cal L}[x^{n-1}R_{n-1}(x)] =
\l_n{\cal L}[x^{n-2}R_{n-2}]$. Taking into account 
the initial condition 
${\cal L}[1] = \l_1$ we establish (2.8) and the proof of the
corollary
is complete.

Note that when $a_n = 0$ for $n \ge 2$ we 
have a Favard theorem associated
with general $T$-fractions with all the moments ${\cal L}[x^n]$,
$n= 0, \pm1, \pm 2,\dots$ uniquely determined \cite{He:Va} 
\cite{Jo:Nj:Th}.

Our next result is a representation theorem 
for the functional $\cal L$ based on the properties of a
continued 
fraction represented as an integral transform. 
Recall that the polynomials  of the second 
kind associated with the recurrence relation (2.1) are
 given by
 \be
 Q_n(x) - (x - c_n) Q_{n-1}(x) + \l_n(x - a_n) Q_{n-2}(x) = 0, n
\ge 2,
 \ee
 with
 \bea
 Q_0(x)  := 0,\, Q_1(x) :=1 \nonumber
 \eea
The ratio $Q_n(z)/P_n(z)$ is the $n$th convergent of the infinite
continued fraction
\be
R_I(z) = \frac{1}{z - c_1}\begin{array}{c} \\ - \end{array}
\frac{\l_2\,(z - a_2)}{z - c_2}
\begin{array}{c} \\ - \end{array} 
 \frac{\l_3(z-a_3)}{z - c_3} \begin{array}{c} \\ - \end{array}
\cdots
\ee
which we assume becomes a finite fraction  when $z = a_k, k \ge 2$.
 This leads to the following definition.

{\bf Definition}. A continued fraction of the type (2.10) will be
referred to as an $R_I$-fraction provided that it terminates when
$z = a_n, n \ge 2$.

We will henceforth assume that $R_I(z)$ converges to a function 
which vanishes at infinity and whose singularity structure is 
given by a finite number of branch cuts and a denumerable 
number of poles so that
\be
R_I(z):=\int_\Gamma \frac{d\alpha(t)}{z - t},\quad
z\in {\cal C}\setminus
{\cal D}.
\ee
Here the measure $d\alpha(t)$ and the
multiple contour $\Gamma$ are taken
in the following generalized sense.
In order to accomodate a pole at $z_k$ with multiplicity $m_k$ and 
 residue $R_k$ we would include in the right side of (2.11) a term
 $$\frac{1}{2\pi i}\int_{|t - z_k|> \epsilon,\,  |z - t|= \epsilon}
 \; \; \frac{R_k\, dt}{(z-t)\, (t - z_k)^{m_k}},$$
with a suitably chosen positive $\epsilon$.
For a branch cut along a contour $\Gamma_j$ with discontinuity
$\a_j'(t)$ we would include a term $$\frac{1}{2\pi i}
\int_{\Gamma_j}\frac{\a_j'(t)\, dt}{z - t},\quad z \notin
\Gamma_j.$$
We will also assume that the 
domain 
\be
{\cal D} = (\cup_{l = 1}^N\,\Gamma_l)
\cup\{z_k: k = 1, 2, \dots\}^{-}, 
\ee
where $\{\cdot\}^{-}$ denotes a closed set.
In other words we assume that
the
continued fraction converges except at the singular points of 
 the function that it represents. 
 
  From Pincherle's theorem \cite{Jo:Th1} we then have the existence
   of a special minimal
   solution $X_n^{(min)}(z)$ satisfying the same recurrence
    relation as $P_n(z)$, namely
 \be
 X_n^{(min)}(z) - (z - c_n)\, X_{n-1}^{(min)}(z) + \l_n\, (z-a_n)\,
 X_{n-2}^{(min)}(z) = 0, \; \; z \notin {\cal D}
 \ee
 but with the minimality condition
 \bea
 \lim_{n \to \infty} X_n^{(min)}(z)/P_n(z) = 0, \quad z\notin {\cal
D}.
 \nonumber
 \eea
 Hence the additional representation
 \be
 R_I(z) = \frac{X_0^{(min)}(z)}{\l_1\, (z - a_1)\,
X_{-1}^{(min)}(z)}
 = \int_{\Gamma}\, \frac{d\a(t)}{z - t}, \quad z \in {\cal
C}\setminus
 {\cal D}.
 \ee
 We will also assume the normalized asymptotics
\be
R_I(z) = \frac{1}{z} \, + \sum_{n = 2}^{\infty} d_n\, z^{-n},\quad
 |z| \to \infty,
 \ee
 so that $\int_{-\infty}^{\infty} d\a(t) = 1$.

A further technical assumption will be that 
\be
\{a_n: 2 \le n \} \cap {\cal D} = \emptyset.
\ee
 From the  viewpoint of multiple point rational interpolants this
means that the sequence of interpolation points associated
 with $R_I(z)$ is $\{\infty, a_2, \infty, a_3, \infty, \dots\}$
  with the $a_n$'s distinct from the singular points of 
   ${\cal D}$, \cite{St:To}.

\noindent{\bf Remark.} Note that $R_I(z)$ in (2.10) does not depend
on $\lambda_1,a_1$ while they occur in (2.14). This seeming
dependence
on $\lambda_1,a_1$ can be interpreted in two ways. Firstly the 
denominator $\lambda_1(z-a_1)X_{-1}^{(min)}(z)$ is determined
through
the recurrence relation (2.13). That is
\bea
\lambda_1(z-a_1)X_{-1}^{(min)}(z) :=
(z-c_1)X_0^{(min)}(z)-X_1^{(min)}(z). \nonumber
\eea
Secondly, in explicit models
where the sequences $\{\lambda_n\}, \{a_n\}$ and
$\{X_n^{(min)}(z)\}$
have an explicit analytic dependence on $n$, there are natural
choices
for $\lambda_1, a_1$ and $X_{-1}^{(min)}(z)$. It often turns out
that,
for these natural values, $\lambda_1=0$. In these cases $(z-a_n)
X_{n-2}^{(min)}(z)$ is singular at $n=1$ and the product $\lambda_1
(z-a_1)X_{-1}^{(min)}(z)$ is indeterminant. For these cases it is
convenient and correct to define
$\lambda_1(z-a_1)X_{-1}^{(min)}(z)$ through
\bea
\lambda_1(z-a_1)X_{-1}^{(min)}(z)
:=\lim_{n\to 1}\lambda_n(z-a_n)X_{n-2}^{(min)}(z). \nonumber
\eea
 This is also true
for (2.18). A similar situation occurs for $R_{II}$-type fractions
in (3.16) and (3.17).

   We now state and prove a representation theorem.
   \begin{thm}
Consider the three term recurrence relation (2.1) and assume
 the representation (2.11) together with conditions (2.12), (2.15)
  and (2.16).  Then the linear functional of Theorem 2.1 has the 
   representation
   \be
   {\cal L}[f] = \int_{\Gamma} \, f(t)\, d\a(t).
   \ee
   \end{thm}
   {\bf Proof}. Let $\l_1 = 1$. Then $\int_\Gamma d\a(t) = 1 = 
    \l_1$. It remains to prove that $\int_\Gamma t^k\, R_n(t)
d\a(t) 
     = 0, \; 0 \le k < n$. This follows from the Lemma below by
taking
      the $z \to \infty$ asymptotics which yields for $0 \le k
       < n$
       \bea
       z^{-1} \int_{\Gamma} t^k\, R_n(t)\, d\a(t) =
       z^{k-n-1}\l_1\, \l_2\cdots \l_{n+1}\, [1 + O(1/z)],
       \nonumber
       \eea
   where we have used the fact that $X_n^{(min)}(z)/[\l_1 \,
   X_{-1}^{(min)}(z)]
   \approx \prod_{j=2}^{n+1}\l_j\; as \,  z \to \infty$ which
follows
    from (2.13)-(2.15). 
    \begin{lem}. Consider the three term recurrence relation
     (2.1) and assume that the representation (2.11) together with
      the conditions (2.12), (2.15) and (2.16) hold. Then we have
   \be
   \frac{z^kX_n^{(min)}(z)}{\l_1\, X_{-1}^{(min)}(z)\,
\prod_{j=1}^{n+1}
    (z - a_j)} = \int_\Gamma \frac{t^kR_n(t)\, d\a(t)}{z - t},
    \; 0 \le k \le n.
    \ee
    \end{lem}
    {\bf Proof}. The singularities in (2.18) come from the zeros 
     of $\lambda_1(z-a_1)X_{-1}^{(min)}(z)$ and the discontinuities
in
     $X_n^{(min)}(z)/\lambda_1(z-a_1)X_{-1}^{(min)}(z)$ across the
contours $\G_j$, that 
    is, the same singularities as in (2.14). The singularities 
    which would seem to appear from the denominator factor 
    $\prod_{j=1}^{n+1} (z - a_j)$ are, without loss of generality,
    canceled by corresponding zeros of $X_n^{(min)}(z)$ (or
    alternatively poles of $\lambda_1X_{-1}^{(min)}(z)$ if one
chooses a different
     normalization). In order to justify the right side of (2.18)
  it remains to compute the "weight" of the singularity. If
 $\lambda_1X_{-1}^{(min)}(z_j) = 0, X_0^{(min)}(z_j) \ne 0$ then it
follows from
 (2.12) and (2.1) that $X_n^{(min)}(z_j) = P_n(z_j)\,
 X_0^{(min)}(z_j)$ since each side of 
 this equality satisfies the same recurrence relation
  and the same initial conditions. Thus for a pole at $z_j$ of
  multiplicity $m_j$ in (2.18) one has a residue 
   $$m_j!\, z_j^k\, P_n(z_j)X_0^{(min)}(z_j)/[\l_1 \prod_{l = 1}^{n
+ 1}
    (z_l - a_l)\left (\frac{d^{m_j}}{dz^{m_j}}
X_{-1}^{(min)}\right)
    (z_j)]$$
    as compared with the residue
    $$X_0^{(min)}(z_j) m_j! /[\l_1\, (z_j - a_1)\left(
    \frac{d^{m_j}}{dz^{m_j}}
    X_{-1}^{(min)}\right)(z_j)]$$
    in (2.14). Hence we see that the pole 
    singularities in (2.18) have a residue
  with an additional factor $z_j^{k} R_n(z_j)$ as required. The
  situation is similar for the absolutely continuous contribution.
  Thus if $\Delta
(X_n^{(min)}(z)/\lambda_1(z-a_1)X_{-1}^{(min)}(z))$ is the
  discontinuity across a contour passing through the point $z$
  then $\Delta ( X_n^{(min)}(z)/\lambda_1(z-a_1)X_{-1}^{(min)}(z))
= P_n(z) \Delta (
  X_0^{(min)}(z)/\lambda_1(z-a_1)X^{(min)}_{-1}(z))$, since each
side of this equality
   again satisfies the same three term recurrence and the same 
   initial conditions at $n = -1, 0$. The condition $0 \le k \le n$
    is required so that (2.18) has at least $O(1/z)$ asymptotics
    and thus a zero at infinity. 

    We now illustrate the above theory with an example which is the
    $q$-analog of Jacobi-Laurent polynomials \cite{He:Va}. 
    Although our example involves a special type of $R_I$-fraction,
  namely a general $T$-fraction, we consider it  to be of intrinsic
  interest.

\bigskip

  {\bf Example 2.1: $_2\phi_1$ Functions}. 
  From the contiguous relations for $_2\phi_1$
   hypergeometric functions in \cite{Is:Li}, or by expanding and
equating
   powers of $z$, it can be verified that
   \be 
   X_{n+1}(z) - \left( z + q^{1/2} \frac{(1 - bq^n)}{(1 -
aq^{n+1})}
   \right) X_n(z) + q^{1/2}z \frac{(1 - q^n)\,(1 -
abq^n)}{(1-aq^{n+1})
   (1 - aq^n)} X_{n-1}(z) = 0
   \ee
   has solutions
   \be
  X_n^{(1)}(z):= z^n
   \frac{(aq^{n+1}, bq^{n+1};q)_{\infty}}{(q^{n+1},
abq^{n+1};q)_\infty}
\,_2\phi_1(1/a, q^{n+1}; b q^{n+1}; q, azq^{1/2})
\ee
and
 \be
 X_n^{(2)}(z):= q^{n/2}\, \frac{(aq^{n+2}, aq^{n+1};q)_\infty}
 {(q^{n+1}, abq^{n+1}; q)_\infty} \, _2\phi_1(q/b, q^{n+1};
aq^{n+2};
  q, bq^{1/2}/z)
 \ee
 and a polynomial solution
 \be
 P_n(z) := \frac {q^{n/2}\, (b;q)_n}{(aq;q)_n} \, _2\phi_1(q^{-n},
 aq; q^{1-n}/b; q, zq^{1/2}/b).
 \ee
 The large $n$ asymptotics of these solutions is easily seen to be
 given by
 \be
 X_n^{(1)}(z) \approx z^n\, (zq^{1/2}; q)_\infty/(az
 q^{1/2};q)_{\infty},
 \ee
\be
X_n^{(2)}(z) \approx q^{n/2}
(q^{3/2}/z;q)_{\infty}/(bq^{1/2}/z;q)_\infty,
 \ee
and
 \be
 P_n(z) \approx \left\{ \begin{array} {c} q^{n/2} \, \frac{(b,
azq^{1/2}; q)_\infty}
  {(aq, zq^{-1/2}; q)_{\infty}},\quad |z| < |q|^{1/2}
  \\ \\ z^n \,
\frac{(bq^{1/2}/z,qb;q)_{\infty}}{(q/b,q^{1/2}/z;q)_{\infty}}
  ,\quad |z| > |q|^{1/2}. \end{array} \right .
  \ee
  Thus the minimal solution to the recurrence relation (2.19)
  is given by
  \be
  X_n^{(min)}(z) = \left \{ \begin{array} {c} X_n^{(1)}(z),
  \quad |z| < |q|^{1/2}
  \\ X_n^{(2)}(z), \quad |z| > |q|^{1/2} .\end{array} \right .
  \ee
  Pincherle's theorem \cite{Jo:Th1} then establishes 
  the continued fraction representation
  \be
  R_I(z) = \left\{  \begin{array} {c}
q^{-1/2}\frac{(1-aq)}{(1-b)}\,
  _2\phi_1(1/a, q; qb; q, zaq^{1/2}),\;\; |z| < |q|^{1/2} \\ \\
  z^{-1}\, _2\phi_1(q/b, q; aq^2; q, bq^{1/2}/z),\; \; |z| >
|q|^{1/2}
  \end{array} \right.
  \ee
 where
 \be  
R_I(z):= \left[z - c_1 - {\bf K}_{k=2}^{\infty}\{\l_k\, z/(z -
c_k)\}
\right]^{-1},
\ee
with
\bea
c_n = -q^{1/2}\frac {(1 - bq^{n-1})}{(1 - aq^n)}, \quad
\quad \l_n =q^{1/2}\frac {(1-q^{n-1})(1 - abq^{n-1})}{(1-aq^n)(1
-aq^{n-1})}.
\eea
Recall the Heine transformation \cite[(III.1)]{Ga:Ra}
\be
_2\phi_1(a, b; c; q, z) = \frac{(b, az; q)_\infty}{(c, z;
q)_\infty}
\; _2\phi_1(c/b, z; az; q, b).
\ee
  From the above transformation we see that the  $_2\phi_1$'s 
  on the right side of (2.27) have singularities at
 $z = q^{-n-1/2}/a$ and $z = b q^{n+1/2}$, $n = 0, 1, \cdots$.
However,
 because of the respective conditions $|z| < |q|^{1/2},\; |z| > 
 |q|^{1/2}$ the right side of (2.27) has no pole singularities if 
 $|aq| < 1$ and $|b| < 1$.  Thus when $|aq| < 1$ and $|b| < 1$ we 
 also have the representation
 \be
 R_I(z) = \int_{|t| = |q|^{1/2}} \frac {\a'(t)\, dt}{z - t}, 
 \quad |z| \ne q^{1/2},
 \ee
 with $\a'$ given by
 \bea
 \a'(t) := \lim_{\epsilon \to 0^+, n \to 0} \frac {1}{2\pi \, i
 \; \l_{n+1}t} \left( \frac {X_n^{(min)}(t_+)}{X_{n-1}^{(min)}
 (t_+)} - \frac{ X_{n}^{(min)}(t_- )}{X_{n-1}^{(min)}
 (t_-)} \right), \, t=|q|^{1/2}e^{i\theta},
\eea
\bea
\quad t_{\pm} := (|q|^{1/2} \pm \epsilon) e^{i\theta}, 0 \le \theta
\le 2\pi. \nonumber
\eea
  A calculation using (2.26),(2.30) and the $q$-Vandermonde
nonterminating
   sum \cite[II.23]{Ga:Ra} yields
 \be
\a'(t) = \frac{i}{2\pi q^{1/2}}\; \frac {(q^{1/2}t,
q^{1/2}/t, q, abq; q)_{\infty}}{(aq^{1/2}t,
bq^{1/2}/t, aq^2, b;q)_\infty}.
\ee
 From our general theory associated with Theorem 2.3 
or the general theory of 
$T$-fractions \cite{Jo:Nj:Th} we may then state the following
results
\be
\int_{|t| = |q|^{1/2}} t^{-k} P_m (t) \, \a'(t)\,  dt = 0, \quad
0 < k \le
m,
\ee
\be
\int_{|t| = |q|^{1/2}} P_m(t)\; \a'(t)\, dt = \frac{q^{m/2}(q,
abq;q)_m}
{(aq, aq^2; q)_m},
 \ee
 and
 \be
- \int_{|t| = |q|^{1/2}} t^{-m-1}\, P_m(t)\, \a'(t)\, dt =
 \frac{(1-aq)q^{-1/2}\, (abq, q; q)_m}{(1-b)\, (aq, bq; q)_m}.
 \ee
 We now derive a rational biorthogonality using (2.34) and (2.36).
    As a first step we shift the contour of integration in (2.34) 
 and (2.36) to the unit circle $|t| = 1$ in order to obtain
 \be
  \int_0^{2\pi} t^{-k}\, P_m(t)\, f(a,b,t)\,
  d\theta = 
  \frac{(q, abq;q)_m}{(aq, bq; q)_m}\delta_{k, m}, 
  \quad 0 \le k \le m, \quad t = e^{i\theta} 
  \ee
  with $f(a, b, t)$  given by
  \bea
  f(a, b, t) = \frac{1}{2\pi}\; \frac{(q^{1/2}t, q^{1/2}/t, q, qab
  ; q)_\infty}{(aq^{1/2}t, bq^{1/2}/t, qa, qb;q)_\infty}. \nonumber
  \eea
  
  The condition for shifting the contour is  that no poles of
   $\a'(t)$ are crossed as  one goes from $|t| =|q|^{1/2}$ to
   $|t| =1$. Thus one needs $|q^{-n+1/2}| > |a|$ and $|bq^{n}| <
1$,
    $n = 1, 2, \cdots$, i.e. $|a| < |q|^{-1/2}$ and $|b| < 1$. The 
    latter two conditions are satisfied since we assumed that $|a|
    < |q|^{-1}$ and $|b| < 1$.

    We now take the complex conjugate of (2.37) followed
    by the replacements $(q, a, b) \to
    (\overline{q}, \overline{b}, \overline{a})$
    to obtain 
\be
 \int_0^{2\pi} t^{k}\, Q_m(1/t) \, f(a, b, t)
  \, d\theta = \frac{(q, abq;q)_m}{(aq, bq;q)_m}\delta_{k,m},\quad
  0 \le k \le m,
  \ee 
  where $Q_m(z)$ is the polynomial
  \be
  Q_m(z) = \frac{q^{m/2}(a; q)_m}{(bq;q)_m}\; _2\phi_1(q^{-m},
  bq; q^{1-m}/a; q, zq^{1/2}/a).
  \ee
  From (2.37) and (2.38) we finally obtain the biorthogonality 
  relation 
  \be
\int_0^{2\pi} P_m(t)\,Q_n(1/t)\,f(a, b, t) \, d\theta =
\frac{(q, abq;q)_m}{(aq, bq;q)_m}\, \delta_{m,n}, \quad |a| < 1, 
|b| < 1.
\ee

This biorthogonality had been previously obtained by Pastro
\cite{Pa}
 using other methods. In \cite{Al:Is} the biorthogonality (2.40)
was
 used to derive a biorthogonality for $_4\phi_3$ rational
functions.
  Here we have shown that the biorthogonality (2.40) is 
  a byproduct of the theory of $R_I$ fractions and we have derived
  a more general transform given by (2.18); namely if $|aq| <1$, 
  $|b| < 1$, $0 \le k \le n$ then 
  \be
  \frac{i}{2\pi}\, \int_{|t| = |q|^{1/2}} \frac{t^{k-n}P_n(t)}
  {z - t}\; \frac{(q^{1/2}t, q^{1/2}/t; q)_\infty}{(aq^{1/2}t, 
  bq^{1/2}/t; q)_\infty} \, dt
  \ee
  \bea
  = \left\{ \begin{array} {c} z^k\frac{(aq^{n+1}, bq^{n+1};
q)_\infty}
  {(q^{n+1}, abq^{n+1}; q)_\infty} \, _2\phi_1(1/a, q^{n+1};
bq^{n+1}
  ; q, azq^{1/2}),\; |z| < |q|^{1/2}  \\  \\
  q^{(n+1)/2} z^{k-n-1} \frac{(aq^{n+1}, aq^{n+2}, b; q)_\infty}
  {(q^{n+1}, abq^{n+1}, aq; q)_\infty}\, _2\phi_1(q/b, q^{n+1}; 
  aq^{n+2}; q, bq^{1/2}/z), \; |z| > |q|^{1/2} \end{array}
  \right. 
  \nonumber
  \eea
  where $P_n(z)$ is as in (2.22). 

  The case $z = k = n = 0$ in (2.41) is equivalent to the original
  $q$-beta integral of Ramanujan with which Pastro
  started \cite{Pa}.

 \bigskip

 \setcounter{section}{3}
 \setcounter{equation}{0}
 
{\bf 3.$R_{II}$-Fractions}.
The format of Section 2 is repeated with first a Favard type
theorem,
 then its realization in terms of the properties of an $R_{II}$ 
 type continued fraction. We include an example of a new system of 
 biorthogonal rational functions. In Section 4 we  cast the recent
results of
\cite {Al:Is} and \cite {Is:Ma} in the language of $R_{II}$
fractions and 
use the theory outlined in this work to give  new derivations
 of these results. We also discover a new set of biorthogonal
rational functions in Section 5. 
 
 Consider the system of polynomials $\{P_n(x)\}$ satisfying
 the recurrence relation
 \be
   P_n(x) - (x - c_n) \, P_{n-1}(x) + \l_n \, (x - a_n)\, 
    (x-b_n)P_{n-2}(x) = 0, \quad n \ge 1,
    \ee
 and the additional assumptions
 \bea
    P_{-1}(x) = 0,\quad  P_0(x) = 1,  \quad \l_{n+1} \ne 0, 
    \quad P_n(a_{n+1}) \ne 0,
    \quad P_n(b_{n+1}) \ne
    0,\quad  n > 0. \nonumber
    \eea
    The recurrence relation (3.1) can be renormalized to yield the
    rational function recurrence relation
    \be
    (x -a_{n+1})\, (x-b_{n+1})\, S_n(x) - (x-c_n)S_{n-1}(x)
    \, + \l_n\,S_{n-2}(x) = 0, n > 1,
    \ee
    with
    \be
    S_n(x) = P_n(x)/\prod_{k = 1}^n [(x - a_{k+1})(x - b_{k+1}).
    \ee
    We now come to an analog of Favard's theorem.
    \begin{thm}
    Given the recursion (3.1) there is a linear functional ${\cal 
    L}$ defined on the span of the rational functions $\{x^k 
    S_n(x)\, : 0 \le k \le n < \infty\}$,  mapping it into
    ${\cal C}$ and normalized by ${\cal L}[1] = N_0$, 
    ${\cal L}[xS_1(x)]
    = N_1$ such that the orthogonality relation ${\cal L}[x^k\,
S_n(x)]
    = 0, 0 \le  k < n$ holds. Furthermore the values of ${\cal L}
    [\prod_{j =1}^n(x-a_{j+1})^{-1}\; \prod_{k=1}^m(x -
b_{j+1})^{-1}]$
    for $m,n = 0, 1, \cdots$, are uniquely determined.   
  \end{thm}
    {\bf Proof}.  Define a linear functional ${\cal L}$ by
requiring 
    it to satisfy
    ${\cal L}[1] = N_0$, ${\cal L}[xS_1(x)] = N_1$, ${\cal
L}[S_1(x)]=0$
    and ${\cal L}[S_n(x)] = {\cal L}[xS_n(x)] = 0$, $n\ge2$. The
recurrence 
    relation (3.2) then yields ${\cal L}[x^k S_n(x)] = 0$, $1 < k
< n$.
    From ${\cal L}[xS_1(x)] = N_1$ and ${\cal L}[S_1(x)]
    = 0$ we obtain ${\cal L}[(x - c_1)/(x - a_2)] =
    {\cal L}[(x-c_1)/(x - b_2)] = N_1$ and hence
    \bea
  {\cal L}[1 +(a_2 - c_1)/(x - a_2)] = {\cal L}[1 + (b_2-c_1)/(x -
b_2)]
  = N_1. \nonumber
  \eea
  With ${\cal L}[1] = N_0$ this implies ${\cal L}[(x - a_2)^{-1}]
=
  (N_1 - N_0)/(a_2-c_1)$ and ${\cal L}[(x-b_2)^{-1}] =
  (N_1-N_0)/(b_2-c_1)$.  If $a_2 \ne b_2$ then a partial fraction 
  decomposition yields $${\cal L}[(x - a_2)^{-1}(x - b_2)^{-1}]
  = (a_2 -b_2)^{-1}\,{\cal L}[(x - a_2)^{-1}-(x-b_2)^{-1}] = (N_0
- N_1)/[(a_2 - c_1)
  (b_2 - c_1)].$$ On the other hand if $a_2 = b_2$
  then ${\cal L}[S_1(x)] = 0$ yields ${\cal L}[(x - a_2)^{-1}
  + (a_2 - c_1)(x - a_2)^{-2}] = 0$, which implies
  ${\cal L}[(x-a_2)^{-2}] = (a_2 - c_1)^{-2}(N_0 - N_1)$. 
  We continue this argument by induction on $n$ using 
${\cal L}[x^kS_n(x)] = 0,\; 0 \le k < n, \; n \ge 2$. Thus for each
new $n, n\ge 2$ we use ${\cal L}[S_n(x)] = {\cal L}[xS_n(x)] = 0$
to 
 evaluate ${\cal L}[\prod_{j=1}^N (x- a_{j+1})^{-1}\, \prod_{k=1}^M
 (x - b_{j+1})^{-1}]$, $N = n$, $M < n$, $N < n$, $M = n$ 
 and $N = M = n$.

The next result gives a recursive definition of ${\cal
L}[x^n\,S_n(x)]$.
\begin{cor}
  Set ${\cal L}[x^n S_n(x)] = N_n$. Then the $N_n$'s satisfy the 
  three term recurrence relation
 \end{cor}
 \be
N_n - N_{n-1} + \l_n\, N_{n-2} = 0, \quad n \ge 2.
\ee
{\bf Proof}. Multiply the recursion (3.2) by $x^{n-2}$, 
apply ${\cal L}$ and make use of
the orthogonality relation ${\cal L}[x^k\,
S_n(x)] = 0, \; 0 \le k < n$. The result is (3.4).
 
 It is because of (3.4) that we need two initial normalizations
 $N_0$ and $N_1$. In Example 3.2, which comes later in this
section,
 we will see that $N_n = \kp_1 \, \kp_2\cdots \kp_{n+1}$
 with $\kp_1 = 1/(1- \kp_2)$. Now $\kp_j = \l_j/(1 - \kp_{j+1}),\;
j \ge 2$ 
 follows from (3.4).  Since the continued fraction ${\bf K}_{n =
2}^{\infty}
 \left( \frac{\l_n}{1}\right)$ of the example converges  we find
that 
 (3.4) can be realized in practice by having
 \be
 \kp_j = K_{n = j}^{\infty} \left( \frac{\l_n}{1} \right),\;
 j > 1, \quad 
 \kp_1 = (1 - \kp_2)^{-1},
 \ee
 and
 \be
 N_n = \kp_1\,\kp_2\, \cdots\, \kp_{n+1}.
 \ee
 
 We now come to a representation theorem for $R_{II}$-fractions.
 The polynomials of the second kind associated with the recursion
 (3.1) are generated by
 \be
  Q_n(x) - (x - c_n)\, Q_{n-1}(x) + \l_n\, (x - a_n)\,(x-b_n)\, 
  Q_{n-2}(x) 
  = 0, \; n > 1,
  \ee
  \bea
  Q_0(x) := 0, Q_1(x) := 1. \nonumber
  \eea
  The ratio $Q_n(z)/P_n(z)$ is the $n$th convergent (approximant) 
  of the continued fraction
  \bea
  R_{II}(z) = \frac{1}{z- c_1} \begin{array} {c} \\  - \end{array}
  \frac{\l_2\, (z - a_2)\, (z - b_2)}{z - c_2} \begin{array} {c}
   \\  - \end{array} 
   \frac{\l_3\, (z - a_3)\, (z- b_3)}{z-c_3} \begin{array}
   {c} \\ -\; \cdots  \end{array},
  \eea
  which we assume becomes a finite fraction for $z = a_k$ or $z =
b_k$
  for any $k >1$. This leads to the following definition.
  
  {\bf Definition}. A continued fraction of the type (3.8) will
  be referred to as  an $R_{II}$-fraction provided that it
  terminates in the cases $z = a_n$ and $z = b_n$, $n \ge 1 $ .

  We will henceforth assume that $R_{II}(z)$ converges to a
function 
  which vanishes at infinity and whose singularity structure is 
  given by at most a finite number of branch cuts $\{\G_j\}_1^N$
and a
  denumerable number of poles $\{z_j\}_1^{\infty}$ so that, in the
  generalized sense previously explained in Section 2,  the
  following representation holds.
  \be
  R_{II}(z) = \int_{\G} \frac{d\a(t)}{z - t}, \quad z \in
  {\cal C}\setminus {\cal D}.
  \ee
  Also, as before, we assume that ${\cal D}$ is minimal so that
  \be
  {\cal D} = \left(\bigcup_{j = 1}^N \G_j\right) \bigcup \{z_k\}^-
  \ee 
  and that the set of interpolation points
  $\{a_n, b_n\}_{n=2}^{\infty}$, see \cite{St:To}, is disjoint from
the
   set of singular points ${\cal D}$, that is
   \be
   \{a_n, b_n\}_{n = 2}^{\infty} \cap {\cal D} = \emptyset.
   \ee
   We will also assume the large $z$ asymptotics
   \be
   R_{II}(z) \approx \kp_1/z = z^{-1}\int_{\g} d\a(t),
   \ee
   with
   \bea
   \kp_1 = \frac{1}{1} \begin{array} {c} \\ - \end{array}
   \frac{\l_2}{1} \begin{array} {c} \\ - \end{array} \frac{\l_3}
   {1} \begin{array}{c} \\ - \cdots  \end{array}
   \eea
   a convergent continued fraction.
   \begin{thm}
   Consider the three term recurrence relation (3.1) and assume
that
   the representation (3.9) and the conditions (3.10)-(3.13) hold.
    Then the linear functional ${\cal L}$ of Theorem 3.1 has the 
   integral representation
  \be
  {\cal L}[f] = \int_{\G} f(t)\, d\a(t).
  \ee
  Furthermore the normalization constants of Corollary 3.2 are
 realized by (3.5) and (3.6).
 \end{thm}
 {\bf Proof}. We set $N_0 = \kp_1$. it remains to show that $N_1
  = \kp_1 - 1$ and that
  \be
  \int_{\G} t^k\, S_n(t) \, d\a(t) = 0, \quad 0 \le k < n.
  \ee
  This will follow from the large $z$ asymptotics in (3.16) of the
  Lemma below.
  \begin{lem}
  Under the assumptions of Theorem 3.3 we have the 
  integral representation
  \be
  \frac{z^k\; X_n^{(min)}(z)}{\l_1\, \prod_{j = 1}^{n+1} [(z - a_j)
  \,(z - b_j)] \; X_{-1}^{(min)}(z)} = \int_\G \frac{t^k\,
S_n(t)}{z -t}
  d\a(t), \; 0 \le k \le n, \; z \in {\cal C}\setminus {\cal D}.
  \ee
  where
  \bea
  \lambda_1(z-a_1)(z-b_1)X_{-1}^{(min)}(z)=(z-c_1)X_0^{(min)}(z)-
  X_1^{(min)}(z) \nonumber
  \eea
  (see the remark before Theorem 2.3).
  In (3.16) $X_n^{(min)}(z)$ denotes the minimal solution to
  the recurrence relation (3.1).
 \end{lem}
 {\bf Proof}. Here again we invoke Pincherle theorem and establish
 the representation
 \be
 \frac{X_0^{(min)}(z)}{\l_1\, (z - a_1)\,(z - b_1) \,
X_{-1}^{(min)}(z)}
 = \int_{\G} \frac{d\a(t)}{z - t}, \quad z \in {\cal C}\setminus
 {\cal D}.
 \ee
The asymptotic relationship (3.12)  gives 
\be
\frac{X_0^{(min)}(z)}{\l_1\,X_{-1}^{(min)}(z)}\approx \kp_1\,z.
\ee

\noindent The three term recurrence relation then implies
\bea 
X_1^{(min)}(z)/X_0^{(min)}(z) \approx (1-1/\kp_1)\, z =\kp_2\, z,
\nonumber
\eea
say, with $\kp_2 := (1 - 1/\kp_1)$ and 
\bea
X_n^{(min)}(z)/X_{n-1}^{(min)}(z) \approx z(1 - \l_n/\kp_n) =
\kp_{n+1}z\;
n \ge 2.   \nonumber
\eea
This establishes the large $z$ asymptotics on the left of (3.16) as
\be 
\frac{z^k\, X_n^{(min)}(z)}{\l_1\,
[\prod_{j=1}^{n+1}(z-a_j)(z-b_j)]
\, X_{-1}^{(min)}(z)} \, \approx z^{k-n -1} \, N_n,
\ee
with $N_n = \kp_1\, \kp_2 \cdots \kp_{n+1}$
and the $\kp$'s are given by (3.1). To establish the equality
(3.16)
 we follow the same procedure as in Lemma 2.4. This means that
 the singularities of (3.16)
  and (3.17)  are the same but (3.16) has the additional weight
factor 
  $t^k \, S_n(t),\, k \le n$.
  
 \noindent  Note that the large $z$ asymptotics of (3.16) yields
  \bea
  \int_{\G} t^k\, S_n(t)\, d\a(t) \, \approx z^{k-n}\, N_n, \quad 
  0 \le k \le n. \nonumber
  \eea
  The choice $k = n =1$ yields $N_1 = \kp_1\, \kp_2 = \kp_1 -1$ 
  but the remaining choices give the orthogonality relations
  $\int_\G t^k \, S_n(t)\, d\a(t) = 0, \, 0 \le k < n$ as required
  by Theorem 3.3. 

\bigskip

{\bf Example 3.1:  Chebyshev Polynomials of Type $R_{II}$}.
This case with constant coefficients is given in \cite {Ma2} with
some misprints which are corrected here. Consider
 the recurrence relation (3.1) of the form
\be
X_{n+1}(z)-(z+\sqrt{ab})X_n(z)+{1\over
4}(z-a)(z-b)X_{n-1}(z)=0,\quad a,b>0.
\ee
This difference equation has solutions 
\be
X_n^{\pm}(z)=\left({(\sqrt{z} \pm \sqrt{a})(\sqrt{z} \pm
\sqrt{b})\over 2}
\right)^n.
\ee

 For $z\notin (-\infty,0]$ the minimal solution is therefore
 \be
X_n^{(min)}(z)=\left({(\sqrt{z}-\sqrt{a})(\sqrt{z}-\sqrt{b})\over
2}\right)^n 
\ee
 with the square root branch chosen so that
 \bea
|z+\sqrt{ab}-(\sqrt{a}+\sqrt{b})\sqrt{z}| <
 |z+\sqrt{ab}+(\sqrt{a}+\sqrt{b})\sqrt{z}|.
\nonumber
\eea
From Pincherle's theorem we then have
\be
\frac{1}{z+\sqrt{ab}}{{}\atop {-}}\frac{(z-a)(z-b)/4}{z+\sqrt{ab}}
{{}\atop {-}}\frac{(z-a)(z-b)/4}{z+\sqrt{ab}}{{}\atop {-}}\cdots =
\frac{2} {(\sqrt{z}+\sqrt{a})(\sqrt{z}+\sqrt{b})}
\ee
\bea
=\frac{2}{ \pi}\int^0_{-\infty}
{(\sqrt{a}+\sqrt{b})\sqrt{-x}dx\over(a-x)(b-x)(z-x)}. \nonumber
\eea
Corresponding to the more general formula (3.16)  we
also have, for $0\le m \le n$,
\be
{P_m(z)X_n^{(min)}(z)(\sqrt{z}-\sqrt{a})(\sqrt{z}-\sqrt{b})
\over (z-a)^{n+1}(z-b)^{n+1}}={1\over \pi}\int_{-\infty}^0
{(\sqrt{a}+\sqrt{b})\sqrt{-x}P_m(x)P_n(x)dx\over
(a-x)^{n+1}(b-x)^{n+1}(z-x)}, 
\ee 
where 
\bea
P_n(z)= \frac{[(\sqrt{z} + \sqrt{a})(\sqrt{z} + \sqrt{b})]^{n+1}
-[(\sqrt{z} -\sqrt{a}) (\sqrt{z}-\sqrt{b})]^{n+1}}
{2^{n+1}\, (\sqrt{a} + \sqrt{b})\sqrt{z}}                               
\eea
 is the polynomial solution
(no longer monic) with initial values $P_{-1}=0,P_0=1$ 

From the large $z$ asymptotics of (3.24) and (3.25) we easily
obtain the orthogonality 
\be
{2\over \pi}\int_{-\infty}^0{(\sqrt{a}+\sqrt{b})P_m(x)P_n(x)
\sqrt{-x}dx
\over (a-x)^{n+1}(b-x)^{n+1}}=2^{-2n+1}(n+1)\delta_{m,n}.
\ee
If we now introduce the rational functions
\be
R_{2n}(x):={P_{2n}(x)\over (a-x)^n(b-x)^n} \quad, \quad 
R_{2n+1}(x):={P_{2n+1}(x)\over (a-x)^n(b-x)^{n+1}}, \quad n \ge 0 
\ee
then (3.24) (with  $m =n$ and $z=a$) and (3.26) (with $n \ne m$)
translate into the rational orthogonality
\be
 {2\over \pi}\int_{-\infty}^0{(\sqrt{a}+\sqrt{b})\sqrt{-x}R_m(x)
R_n(x)dx\over(a-x)^2(b-x)}={2^{-2n}a^{-1/2}\over(\sqrt{a}+
\sqrt{b})}\delta_{m,n}.
\ee

\bigskip
 
 {\bf Example 3.2: $_2F_1$ Functions}.
 A contiguous relation for a hypergeometric function
 of type $_2F_1$ in \cite[(45), p. 104]{Er:Ma1} tells us that
  the recurrence relation
  \be
  X_n(z)= \left(z - \frac{n + a-1}{2n + a -1 -b}\right)\,
X_{n-1}(z)
  \ee
\bea
{} \qquad - \,\frac{z\,(z-1)\, (n-1)\,(n+a-1-b)}{(2n + a -1 -b)\,
  (2n + a -3 -b)} X_{n-2}(z)= 0 \nonumber
  \eea
  has solutions
 \be
 X_n^{(1)}(z) := \left( \frac{z}{2}\right)^n \, \frac{\G(n+1)\, 
 \G(n+1+a-b)}{\G(n+a+1)\,\G(n + (1 + a -b)/2)}\, _2F_1(a, b; n+1
+a;
  z),
  \ee
  and
 \be
 X_n^{(2)}(z) := \left(\frac{z-1}{2}\right)^n \frac{\G(n+1)\,
\G(n+1 +
 a - b)}{\G(n+2-b)\, \G(n+(1+a -b)/2)} \, _2F_1(1-a, 1-b; n+2 -b; 
 1-z)
 \ee
 and a polynomial solution
 \be
 X_n^{(3)}(z) = P_n(z; a, b) := \frac{2^{-n}\,
(1-b)_n}{((1+a-b)/2)_n}
 \; _2F_1(-n, b-a-n; b-n; 1-z).
 \ee

 The asymptotic behavior of $X_n^{(1)}(z)$ and $X_n^{(2)}(z)$ as
 $n \to \infty$ can
 be easily calculated and it yields
 \be
 X_n^{(1)}(z) \approx (z/2)^n n^{-(a+b-1)/2}, \quad |z| < 1,
 \ee
 \be
 X_n^{(2)}(z) \approx \left(\frac{z-1}{2}\right)^n n^{(a + b
-1)/2},
 \quad |z - 1| < 1.
\ee

Thus the minimal solution of (3.29) is given by
\be
X_n^{(min)}(z) = \left\{ \ba{c} X_n^{(1)}(z), \quad {\cal R}z < 1/2
\\ X_n^{(2)}(z), \quad {\cal R}z > 1/2.\ea \right.  
\ee
 Exploiting Pincherle's theorem  we obtain
\bea
R_{II}(z) = \left\{ \ba {c} -\frac{(1 + a -b)\,(1-z)^{b-1}}
{a}\; _2F_1(a, b; 1+a; z), \quad {\cal R}z < 1/2 \\ \\
- \frac{(1+a - b)}{b-1}\, z^{-a} \, _2F_1(1-a, 1-b; 2-b; 1-z),
{\cal R}z 
> 1/2 \ea \right. 
\eea
 where 
 \be
 R_{II}(z) = [ z - c_1 + {\bf K}_{n = 2}^\infty\{z\,(z-1)\,
 \l_n/(z - c_n)\}]^{-1},
 \ee
 \bea
 c_n := \frac{(n+a-1)}{(2n + a -b - 1)}, \quad \l_n := 
 \frac{(n-1)\, (n+a -1-b)}{(2n + a -1 -b)\,(2n + a -3-b)}.
\nonumber
 \eea

 We now have an explicit example of an $R_{II}$ type continued
fraction
 which converges except when ${\cal R}z = 1/2$. To calculate
 the absolutely continuous measure which is now supported on 
 ${\cal R}z =1/2$ we use \cite[2.9(33), p. 107]{Er:Ma1},  
  \be
  \frac{(1+a-b)}{a}(1 - z)^{b-1}\, _2F_1(a, b; 1+a; z)
  - \frac{(1 + a -b)}{(b-1)}\, z^{-a}\, _2F_1(1-a, 1-b; 2-b;
  1-z)
  \ee
  \bea
  \quad \quad = z^{-a}\,
  (1 - z)^{b-1} \frac{(1 + a - b)\,\G(a)\,\G(1-b)}
  {\G(1 + a - b)}.  \nonumber
  \eea
  Thus we have established
  \be
  R_{II}(z) =\frac{(1 + a - b)\, \G(a)\, \G(1-b)}{
  2\pi\, i\, \G(1+ a - b)} \int_{\frac{1}{2} - i\infty}^{
  \frac{1}{2} + i \infty} \frac{t^{-a}\, (1 - t)^{b-1}}{z - t}
  \; dt, {\cal R}z \ne 1/2,
\ee
with the correct $O(1/z)$ asymptotics of the present section when
${\cal R}(a - b) > 0$ which is \cite[2.1.4(17), p. 63]{Er:Ma1}
 \be
 R_{II}(z) \approx \left( \frac{a - b + 1}{a - b}\right)z^{-1}
 , \quad {\cal R}(
 a - b) > 0.
 \ee
This justifies the assumption (3.12) since we may 
rewrite the $\l_n$'s in the form  
 \be
 \l_n = \frac{1}{4} + \frac{1 - (a - b)^2}{4 \,
 (n + (a-b - 3)/2))\, (n + (a - b - 1)/2)}
 \ee
 and a result in \cite{Th:Wa} yields
\bea
\kp_1= \frac{1}{1} \ba {c} \\ - \ea \frac{\l_2}{1} \ba {c} \\ -
\ea \frac{\l_3}{1} \ba{c} \\ - \cdots \ea = \frac{a-b +1 }
{a - b}, \; {\cal R}(a-b) > 0.
\eea
We now derive a biorthogonality relation for certain rational
functions using the orthogonality relation (3.15) with 
  the special case $k \le n, z = 0$ of the integral
 formula (3.16).  We set
 \be
 U_n(x; a, b) := \frac{P_n(x; a, b)}{(x-1)^{n}},\quad
 V_n(x; a. b) := U_n(1-x; -b, -a), 
 \ee
 \be
 g(x; a, b) := x^{-a
 -1}\, (1 - x)^{b-1}/2\pi i.
 \ee
 The orthogonality relation (3.15) now takes the form
 \be
 \int_{\frac{1}{2}-i\infty}^{\frac{1}{2}+i\infty}t^{-m} \, U_n(
 t; a, b)\, g(t; a, b) \, dt = 0, \quad 0 \le m < n.
 \ee
 Taking the complex conjugate of (3.45), replacing $a$ by $-
 \overline{b}$ and $b$ by $-\overline{a}$, taking into account
 that $\overline{t} = 1 - t$ and using (3.16) (with $z^k$ replaced
by
 $P_k(1-z;-b,-a), k=n$ and $z = 0$ )
 we see that
 \be
 \int_{\frac{1}{2} - i\infty}^{\frac{1}{2}+ i\infty} U_n(t; a, b)
 \, V_m(t; a, b)\, g(t; a, b)\ dt = \frac{\G(1 + a -b)\, n!\, (1+
a
 - b)_n}{\G(1+ a)\, \G(1 - b)\, 2^{2n}\, [((a - b +1)/2)_n]^2}\;
 \delta_{m,n},
 \ee
 provided that $a, b \ne 0, a \ne -1, b \ne 1, {\cal R}(a - b) >
0$.

The biorthogonal rational functions in this example differ from the
ones 
considered by  Askey in \cite {As}. Both systems of rational
functions are 
biorthogonal with respect to a Cauchy beta integral weight
function.

\bigskip

\setcounter{section}{4}

\setcounter{equation}{0}

{\bf 4. The ${}_4 \phi_3$ Functions}.
 In this section we provide two examples of biorthogonal rational
functions which have explicit representation as ${}_4\phi_3$'s. The

three term recurrence relations and the continued fractions 
 associated with these systems come from  contiguous relations for
  ${}_4\phi_3$ and $_8\phi_7$ function. The  ${}_4\phi_3$
contiguous relation is the same one used
 by Askey and Wilson \cite{As:Wi} and that led to the Askey-Wilson
 polynomials. The  difference between our examples and the
 Askey-Wilson polynomials is a reparametrization which converts the
three term recurrence to one of type $R_{II}$. 

 The monic Askey-Wilson polynomials are defined by \cite{As:Wi}
 \bea
 P_n(x; \a, \b, \g, \delta |q) := 
 \frac{(\a\b, \a\gamma , \a \delta, \a\b\g \delta /q;q)_n}
 {(2\a)^n \; (\a\b\g \delta/q; q)_{2n}}
 \;{}_4\phi_3\left( \left. \ba{c}
 q^{-n}, \a\b\g\delta q^{n-1}, \a u, \a /u \\
 \a\b, \a \g, \a \delta \ea \right| q, q\right),
 \eea
 where 
 \bea
 x := \frac{1}{2}(u + 1/u).
 \eea
 The $P_n$'s satisfy the three term recurrence relation
 \be
P_n(x) -(x - a_{n-1})\, P_{n-1}(x) + A_{n-2}\, B_{n-1}\, P_{n-2}(x)
= 0,
 \ee
 \bea
 a_n = \frac{1}{2}(\a + 1/\a) \, - A_n\, - B_n, \nonumber
 \eea
 \bea
 A_n : =
 \frac{(1 - \a\b\g\delta q^{n-1})(1 - \a\b q^n)(1-\a\g q^n)
 (1 - \a\delta q^n)}
 {2\a \; (1 - \a\b\g\delta q^{2n-1})\, (1 -
 \a \b \g \delta q^{2n})},
 \nonumber
 \eea 
 \bea
 B_n = \frac{\a (1-q^n)(1 - \b \g q^{n-1})(1 - \b\delta q^{n-1})(1
- 
 \g\delta q^{n-1})}
 {2\; (1 - \a\b\g\delta q^{2n-2})\;(1 -\a\b\g\delta
   q^{2n-1})}.  \nonumber
  \eea
 
 In \cite{Gu:Ma} the minimal solution to the recurrence relation 
  (4.3) for $n > 0$ has been shown to be
  \bea
  Y_n^{(min)}(x) = \left\{ \ba{c} F_n(u), \quad |u| > 1,
  \\ F_n(1/u), \quad |u| < 1,
   \ea  \right.
   \eea
  and the $F_n$'s are the functions
  \be
  F_n(u) := \frac{(2u)^{-n}\, (\a\b\g\delta q^{2n-1}, \a q^{n+1}/u,
  \b q^{n+1}/u, \g q^{n+1}/u, \delta q^{n+1}/u; q)_\infty}
  {(q^{n+1}, q^{n+2}/u^2, \a\b q^n, \a\g q^n, \a\delta q^{n}, \b\g
q^n,
   \b\delta q^n, \g \delta q^n; q)_\infty}
   \ee
   \bea
   \quad \quad\times{} _8W_7\left(q^{n+1}/u^2; q^{n+1},
   \frac{q}{\a u}, \frac{q}{\b u},
   \frac{q}{\g u}, \frac{q}{\delta u}; \a\b\g\delta q^{n-1}\right).
   \nonumber
   \eea
  From (4.4) one can calculate, using Pincherle's theorem, the 
  Askey-Wilson weight function
  \be
  w(x) := \frac{1}{2\pi\sqrt{1 - x^2}} 
   \frac{(u^2, 1/u^2, \a \b, \a \g, \a \delta, \b\g, \b\delta, 
   \g\delta,q; q)_\infty}
   {(\a u, \a /u, \b u, \b /u, \g u, \g /u, \delta u, \delta
/u,\a\b\g\delta;
   q)_\infty},
   \ee
 normalized by 
 \bea
 \int_{-1}^{1} w(x)\, dx = 1, \quad |\a |< 1, |\b| < 1, |\g| < 1, 
  |\delta| < 1. \nonumber
  \eea
   The system of biorthogonal rational functions recently given by
    Al-Salam and Ismail \cite{Al:Is} can be related to the above
    Askey-Wilson case through a change in parameterization. This
will 
    give us our first explicit $_4\phi_3$ example of an
$R_{II}$-fraction.
  
\bigskip

 {\bf Example 4.1: Biorthogonality on the Unit Circle}. Let us  use
the new parameters $a, b, t_1, t_2$ and the variable $z$
 given by the replacements
 \be
 \a = q^{-1/4}t_2 \sqrt{z}, \; \b = q^{3/4}a \sqrt{z}, \;
 \g =b q^{1/4}/\sqrt{z}, \; \delta = q^{-3/4}t_1/\sqrt{z}, \;
  u = q^{-1/4}\sqrt{z}.
  \ee
  Then after renormalization (4.3) becomes a recurrence relation 
    of type $R_{II}$, namely
    \be
P_n(z) - (z - c_n)\, P_{n-1}(z) + \l_n\, (z - a_n)\,(z- b_n)\, 
P_{n-2}(z)\,=\,0
\ee
with
\be
\l_{n+1} = \frac{at_2q^{n}(q^n-1)(1-abt_1t_2q^{n-2})
 (1 - bt_2q^{n-1})(1 -
 t_1t_2q^{n-2})(1 - abq^n)(1 - at_1q^{n-1})}
 {4\,\sqrt{q} (1 - abt_1t_2q^{2n-3})\,(1 - abt_1t_2q^{2n -2})^2 \,
 (1 - abt_1t_2q^{2n-1})\, u_n\, u_{n+1}}
 \ee
 \bea
a_{n+1} = bt_1q^{n-3/2}, \quad b_{n+1} = q^{-n+1/2}/(at_2),
 \quad c_n = -v_n/u_n, \nonumber
 \eea 
\bea
u_{n+1} = \frac{q^{-1/4}}{2(1 - abt_1t_2q^{2n-2})(1 -
abt_1t_2q^{2n})}
 \nonumber 
 \eea
 \bea
 \times \{1 - q^{n-1}[(1 + \frac{abt_1t_2}{q}q^{2n})
 (qt_2+aq^2+at_1t_2
 + abt_2q) -\frac{abt_1t_2}{q} q^{n}(1 + q)(t_2+qa+\frac{q}{b}+
 \frac{q^2}{t_1})]\}
 \nonumber
 \eea
 \bea 
 v_{n+1} = \frac{q^{1/4}}{2(1-abt_1t_2q^{2n-2})(1 -
abt_1t_2q^{2n})}
 \nonumber
 \eea
 \bea
 \times \{1- q^{n-1}[(1 + \frac{abt_1t_2}{q}q^{2n})(bq + t_1 +
abt_1
 +\frac{bt_1t_2}{q})
 -\frac{abt_1t_2}{q} q^{n}(1 + q)(b + \frac{t_1}{q} + \frac{q}{t_2}
 + \frac{1}{a} )]\}. \nonumber
\eea
 Note that, in order to calculate the middle coefficient in (4.8)
and arrive at values for $c_n$, $u_n$ and $v_n$, it is very useful
to use the form of 
the Askey-Wilson three term recurrence relation given in
 \cite [(1.24), (1.27)]{As:Wi}. This is true also for the next
example.

The polynomials $\{P_n(z)\}$ are now given by
    \bea
    P_n(z)= \frac{q^{n/4}(at_2zq^{1/2}, 
    bt_2, t_1t_2/q, abt_1t_2/q;q)_n}
    {(2t_2)^n(abt_1t_2/q;q)_{2n}\, \prod_{k =1}^n u_k}
    \eea
    \bea
   {}\quad \quad \times{}_4\phi_3\left(\left.\ba{c} q^{-n},
    abt_1t_2q^{n-1}, t_2zq^{-1/2}, t_2
    \\at_2zq^{1/2},\; bt_2,\; t_1t_2/q \ea \right|q,q \right).
   \nonumber
   \eea
   From (4.4) with a renormalization factor
    $z^{n/2}/\prod_{k = -1}^{n}u_k$ we obtain the minimal solution
    to (4.8) 
    \bea
    X_n^{(min)}(z) = \left\{ \ba{c}X_n^{(1)}(z), \quad |z| >
|q|^{1/2}
    \\ X_n^{(2)}(z), \quad |z| < |q|^{1/2}, \ea \right.
    \eea
    with
 \be
 X_n^{(1)}(z) = \frac{(q^{1/4}/2)^n(abt_1t_2q^{2n-1}, t_2q^{n+1},
 aq^{n+2}, bq^{n+3/2}/z, t_1q^{n+1/2}/z;q)_\infty}
 {(q^{n+1}, q^{n+5/2}/z, at_2q^{n+1/2}z, bt_2q^n, t_1t_2q^{n-1},
 abq^{n+1}, at_1q^n, bt_1q^{n-1/2}/z;q)_\infty}
 \ee
 \bea
 \quad \quad \times\,_8W_7\left(\frac{q^{n+3/2}}{z};
  q^{n+1}, \frac{q^{3/2}}{t_2z},
 \frac{q^{1/2}}{az},\frac{q}{b}, \frac{q^2}{t_1};
 abt_1t_2q^{n-1}\right)/\prod_{k = -1}^{n}u_k,  \nonumber
 \eea
 and
 \be
 X_n^{(2)}(z) = \frac{(q^{-1/4}z/2)^n(abt_1t_2q^{2n-1}, 
 t_2zq^{n+1/2}, azq^{n+3/2}, bq^{n+1}, t_1q^n;q)_\infty}
 {(q^{n+1}, zq^{n+3/2}, at_2zq^{n+1/2}, bt_2q^n, t_1t_2q^{n-1},
 abq^{n+1}, at_1q^n, bt_1q^{n-1/2}/z;q)_\infty}
 \ee
 \bea
 \quad \quad \times\,_8W_7\left(q^{n+1/2}z; q^{n+1}, \frac{q}{t_2},
 \frac{1}{a},
 \frac{q^{1/2}z}{b}, \frac{q^{3/2}z}{t_1};
abt_1t_2q^{n-1}\right)/\prod_{k=-1}^{n} u_k.
 \nonumber
 \eea
  With the above minimal solution, equations (3.9) and (3.17) imply
\be
R_{II}(z) =\frac{2u_1q^{-1/4}(1-abt_1t_2/q)(1-q^{1/2}z)}
{(1-b)(1-t_1/q)(1-t_2zq^{-1/2})(1-azq^{1/2})}
\ee
\bea
\quad \quad \times\,_8W_7\left(zq^{1/2}; q,\frac{q}{t_2},
\frac{1}{a},
\frac{q^{1/2}z}{b},\frac{q^{3/2}z}{t_1};abt_1t_2/q\right),\; for\,
|z|
<|q|^{1/2}   \nonumber
\eea
 \be
 R_{II}(z) = \frac{2u_1q^{1/4}(1-abt_1t_2/q)(1-q^{3/2}/z)}
 {z(1-t_2)(1-aq)(1-bq^{1/2}/z)(1-t_1q^{-1/2}/z)}
 \ee
  \bea
 \quad\quad\times\,_8W_7\left(\frac{q^{3/2}}{z}; q,
\frac{q^{3/2}}{t_2z},
  \frac{q^{1/2}}{az}, \frac{q}{b}, \frac{q^2}{t_1};
abt_1t_2/q\right), \; for
  \, |z| > |q|^{1/2}.   \nonumber
  \eea
An application of Cauchy's theorem gives us the integral
 representation
\be
R_{II}(z) = \int_{|t| = |q|^{1/2}}\frac{\a'(t)}{z - t}\, dt,
\ee
with $\a'(t)$ as $\frac{1}{2\pi\,i}$ times the difference of the 
boundary values of $R_{II}(z)$ in (4.11) as $|z| \to |q|^{1/2}$. 
The identity \cite[(III.37), p.246]{Ga:Ra} then gives $\a'(t)$
explicitly. The result is that $\a'(t)$ is $u_1/t^{1/2}$ times the 
Askey-Wilson weight (4.6) but with the new parameterization given
in (4.2) and (4.7), with $z$ replaced by $t$. In order to isolate
the symmetric terms in $\a'$ we chose to write it in the form
\be
\a'(t) = \frac{iu_1}{\pi q^{1/4}}\, \frac{(1-bt_1q^{-1/2}/t)}
{(1-t_1/q)(1-b)}\, f(t)
\ee
where $f(t)$ is
\be
f(t) = \frac{(q^{1/2}t, q^{1/2}/t, at_2q^{1/2}t, bt_1q^{1/2}/t, 
bt_2, at_1, abq, t_1t_2/q,q;q)_\infty}
{(aq^{1/2}t, bq^{1/2}/t, aq, bq, t_1, t_2, q^{-1/2}t_2t,
q^{-1/2}t_1/t, abt_1t_2; q)_\infty}.
\ee
Note that $f(t)$ is symmetric under the transformation
\be
(a, b, t_1, t_2, t)  \; \to \; (b, a, t_2, t_1, 1/t).
\ee

We now derive the biorthogonality relation of Al-Salam and Ismail
\cite{Al:Is}. In the case under consideration our 
orthogonality relation (3.15) becomes
\be
\frac{iu_1}{\pi\,q^{1/4}}\int_{|t| = |q|^{1/2}}
\frac{t^k\,P_n(t)\, (1-bt_1q^{-1/2}/t)\, f(t)\,dt}
{t^n\, (bt_1q^{-1/2}/t;q)_n\,
\prod_{j=1}^n[t - q^{-j+1/2}/(at_2)]} = 0,\; 0 \le k <n.
\ee
We follow the notation in \cite{Al:Is} and define rational
functions
$r_n$ and $s_n$ by
\bea
 r_n(z) = r_n(z; a, b, t_1, t_2) := {}_4\phi_3\left(\left.\ba{c}
q^{-n},
abt_1t_2q^{n-1}, q^{-1/2}t_2z, t_2   \\
at_2q^{1/2}z, bt_2, t_1t_2/q \ea \right| q, q \right),
\eea
and
\be
s_n(z) = s_n(z; a, b, t_1, t_2) := r_n(z; b, a, t_2, t_1).
 \ee
There is a slight difference between the above  definition
of $s_n$ and the definition of $s_n$ in \cite{Al:Is}, due to
complex conjugation. Now the relationship between $P_n$ and $r_n$
is
\be
r_n(z; a, b, t_1, t_2) = \frac{C_n\,P_n(z)}{(at_2q^{1/2}z; q)_n},
\ee
where the normalization constant $C_n$  is given by
\bea
C_n = \frac{(2t_2q^{-1/4})^n\, (abt_1t_2/q;q)_{2n}\prod_{k =1}^n
u_k}{(bt_2, t_1t_2/q, abt_1t_2/q;q)_{n}}. \nonumber
\eea
We may now rewrite (4.20) as
\be
\int_{|t| = |q|^{1/2}} r_n(t)\, s_k(1/t)\, f(t) \, \frac{dt}{t}
= 0, \quad k < n.
\ee
If we deform the contour in (4.23) to the unit circle $|t| = 1$ and
combine this new integral with its complex conjugate (for which 
$\overline{t} = 1/t$) then perform the replacements $(a, b, t_1,
t_2, q) \to (\overline{b}, \overline{a}, \overline{t_2}, 
\overline{t_1}, \overline{q})$ while taking into account the
symmetry
of $f$ under the transformation (4.19) we obtain the
biorthogonality
relation
\be
\int_{|t| = 1}r_k(t; a, b, t_1, t_2)s_n(1/t; a, b, t_1, t_2)
\frac{dt}{t} = 0, \quad k \ne n, \quad n, k \ge 0. 
\ee
The value of this last integral when $n = k$ is deduced 
from a modification of (3.16) with 
 $k = n$ and $z = bt_1q^{n-1/2} < 1$. 
Namely,
 \be
 \frac{z^n\tilde P_n(1/z)X_n^{(min)}(z)}{\l_1
 [\prod_{j=1}^{n+1}(z - a_j)(z - b_j)]X_{-1}^{(min)}(z)} =
\int_{|t| = 
 |q|^{1/2}} \frac{t^n \tilde P_n(1/t)P_n(t) \, \a'(t)\, dt}
 {[\prod_{j = 2}^{n+1} (t - a_j)(t - b_j)] \,(z - t)},
 \ee
 where $\tilde P_n(z) = P_n(x; b, a, t_2, t_1)$.
 We shift the above contour to $|t| = 1$ and put $z =
bt_1q^{n-1/2}$.
 At this value of $z$ the ${}_8W_7$ in $X_n^{(min)}(z)$ is summable

 using \cite[III.23, p. 243]{Ga:Ra}. After a lengthy calculation
 we arrive at 
 \be
 \frac{i}{2\pi} \int_{|t| = 1} r_n(t) s_n(1/t)\, f(t) \,
\frac{dt}{t}
 = -\frac{(t_1t_2/q)^n\, (q, abq, abt_1t_2q^{n-1}; q)_n} 
 {(t_1t_2/q; q)_n(abt_1t_2;q)_{2n}}, 
\ee
in agreement with \cite{Al:Is}.

The condition for the absence of the discrete spectrum and the
validity
of the integrals on the contour $|t| = |q|^{1/2}$ is $|a|, |b|,
|t_1|,
|t_2| < 1$. The condition for shifting the contour to $|t| =1$ is 
to have $|a|, |b| < |q|^{-1/2}$ and $|t_1|, |t_2| < |q|^{1/2}$.
Thus
our overall condition for deriving the biorthogonality is $|a|, |b|
< 1$ and $|t_1|, |t_2| < |q|^{1/2}$.

\bigskip

 {\bf Example 4.2: Biorthogonality on the Line}. We next   use a
different set of parameters $ t_1, t_2, t_3, t_4$ and 
 variable $z$. We set
 \bea
 \a =  \sqrt{t_3t_4} t_1/q, \; \b = -qe^\xi/ \sqrt{t_3t_4}, \;
\g = qe^{-\xi}/ \sqrt{t_3t_4}, \; \delta = t_2q^{-2}/\sqrt{t_3t_4},
\eea
\bea
  u = -\sqrt{t_3/t_4}, \;
  z = \sinh \xi. \nonumber
 \eea
 The explicit representation
 (4.1), after renormalization by $\prod_{k=1}^n u_k$ is now 
\bea
 P_n(z; t_1, t_2, t_3, t_4 |q) &=& 
 \frac{(-t_1e^\xi, t_1e^{-\xi}, t_1t_2t_3t_4q^{-3},
-t_1t_2q^{-2};q)_n}
 {(2t_1\sqrt{t_3t_4}/q)^n \; (-t_1t_2/q^2; q)_{2n} \prod_{k = 1}^n
u_k}\\
 &{}&\times {}_4\phi_3\left( \left. \ba{c}
 q^{-n}, -t_1t_2q^{n-2}, -t_1t_3/q, -t_1t_4/q \\
 -t_1e^{\xi},\;  t_1e^{-\xi}, \; t_1t_2t_3t_4q^{-3} \ea \right| q,
q\right), \nonumber
 \eea
where $u_n$  is given below in (4.33). 

 The three term  recurrence relation (4.3), after renormalization
is now 
    of type $R_{II}$ and is
 \be
P_n(z) - (z - c_n)\, P_{n-1}(z) + \l_n\, (z - a_n)\,(z- b_n)\, 
P_{n-2}(z)\,=\,0
\ee
with
\be
\l_{n+1} = \frac{t_1t_2q^{2n-3}(1- q^n)(1+ t_1t_2q^{n-3})(1 +
q^{n+1}/t_3t_4)
 (1 - t_1t_2t_3t_4q^{n-4})}
 { (1 + t_1t_2q^{2n-2})\,(1  + t_1t_2q^{2n-3})^2 \,(1 +
t_1t_2q^{2n-4})
\, u_n\, u_{n+1}}
 \ee
 \bea
a_{n+1} =\frac{1}{2}(t_2q^{n-2} - q^{2-n}/t_2), \quad b_{n+1} =
\frac{1}{2}
(t_1q^{n-1} - q^{-n-1}/t_1), 
\eea
\bea
  c_{n+1} = \frac{-v_{n+1} + \frac{1}{2}(\sqrt{t_3/t_4} +
\sqrt{t_4/t_3})}{u_{n+1}}, \nonumber
 \eea 
where $u_n$ and $v_n$ are given 
\bea
u_{n+1} &=&
 [(1 - t_1t_2q^{2n-2})(q^4 +t_1t_2t_3t_4) + q^n (1+q)t_1t_2(q -
t_3t_4)]
\\
&& \qquad \qquad\times 
\frac{q^{n-3}/\sqrt{t_3t_4}}{(1 + t_1t_2q^{2n-3})(1 +
t_1t_2q^{2n-1})} \nonumber
 \eea
 \bea 
 v_{n+1} &=& - [(1 - t_1t_2q^{2n-1})(t_3t_4 - q^2) 
 + q^{n-3}(1+q)(t_1t_2t_3t_4 + q^4)] \\
&& \qquad \times  \frac{q^{n-1}(t_1+t_2/q)}{2\sqrt{t_3t_4} (1 +
t_1t_2q^{2n-3}(1 + t_1t_2q^{2n-1})}. \nonumber
\eea

In the case under consideration the continued fraction is the
$R_{II}$ fraction \bea
 \frac{1}{z-c_1}  {{}\atop {-}} \frac{\l_2(z-a_2)(z-b_2)}{z - c_2}
{{}\atop{-}}
\cdots. 
\eea
When the continued fraction converges it will converge to $F(z)$,
where
\bea
F(z) &=& -\frac{2u_1}{\sqrt{t_3/t_4}}\frac{(1-
qt_4/t_3)(1+t_1t_2/q^2)}
{(qt_4/t_3, -t_1t_4/q, qe^\xi/t_3, -qe^{-\xi}/t_3,
-t_2t_4/q^2)_\infty}
 \\
&&
\qquad  \times \tilde{W}(qt_4/t_3; q, -q^2/t_1t_3, t_4e^{-\xi},
-t_4e^{\xi},
-q^3/t_2t_3; t_1t_2/q^2), \nonumber
\eea
with
\bea 
\tilde{W}(a; b, c, d, e, f; \frac{a^2q^2}{bcdef}) := (b, c, d, e,
f;q)_\infty 
{}_8W_7(a; b, c, d, e, f; \frac{a^2q^2}{bcdef}) \nonumber
\eea
The singularities of (4.36) are  at $z = z_n$,
\bea
z_n = \frac{1}{2}(t_3q^{-n-1}-q^{n+1}/t_3), \quad n = 0, 1, \cdots
.\nonumber
\eea
given by the zeros of $(qe^{\xi}/t_3, -qe^{-\xi}/t_3; q)_\infty$.
The residues
 at these points may be calculated explicitly because the
contribution from 
$\tilde{W}$ is in terms of a very well-poised $_6\phi_5$, which is
summable.
  After a straightforward but rather lengthy calculation, we arrive
at the 
Mittag-Leffler expansion
\be
F(z) = \sum_{k = 0}^\infty  \frac{\omega_k}{z - z_k} 
\ee
with
\bea
\omega_k &=& \frac{u_1\sqrt{t_3t_4}}{q^{4k+1}}\;
\frac{(t_1t_2t_3t_4q^{-3},
 qt_1/t_3, 
t_2/t_3, q t_4/t_3; q)_\infty}{(-t_1t_2/q, - t_1t_4/q, -t_2t_4/q^2,
 -q^2/t_3^2; q)_\infty}\\ 
&& \times \frac{(-q^2/t_1t_3, -q^3/t_2t_3, -q^2/t_3t_4, -q^2/t_3^2;
q)_k}
{(qt_1/t_3, t_2/t_3, qt_4/t_3, q; q)_k} (1 +
q^{2k+2}/t_3^2)(t_1t_2t_3t_4)^k. 
\nonumber
\eea
The resulting orthogonality (3.15) is now given by
\be
\sum_{k=0}^\infty \frac{z_k^m\, P_n(z_k)\, \omega_k}
{\prod_{j = 1}^n(z_k- a_{j+1})(z_k- b_{j+1})}= 0, \quad 0 \le m <
n.
\ee

We now derive a rational biorthogonality relation. Let us include
the factor 
\bea
z_k - a_2 = \frac{t_3}{2}q^{-k-1}(1 - t_2q^k/t_3)(1 +
q^{k+2}/t_2t_3)
\nonumber
\eea
in the weight  to obtain the new weight 
\bea
r_k &=& \frac{\omega_k}{z_k - a_2} \ \nonumber \\
&=&\frac{2t_2u_1\sqrt{t_3t_4}}{q^{3k}(q^2+t_2t_3)}\;
\frac{(t_1t_2t_3t_4q^{-3},
 qt_1/t_3, 
qt_2/t_3, q t_4/t_3; q)_\infty}{(-t_1t_2/q, - t_1t_4/q^2,
-t_2t_4/q^2,
 -q^2/t_3^2; q)_\infty} \nonumber \\ 
&&\times \frac{(-q^2/t_1t_3, -q^2/t_2t_3, -q^2/t_3t_4, -q^2/t_3^2;
q)_k}
{(qt_1/t_3, qt_2/t_3, qt_4/t_3, q; q)_k} (1 +
q^{2k+2}/t_3^2)(t_1t_2t_3t_4)^k. 
\nonumber
\eea
Note that this has $k$ dependent terms which are now symmetric in
$t_1$ and $t_2$. This symmetry  together with (4.39) yields the
biorthogonality relation
\be
\sum_{k = 0}^\infty \tilde{R}_m(z_k)\, R_n(z_k)\, r_k = 0, \quad m
\ne n,
\ee
where
\bea
R_n(z) &=& R_n(z; t_1, t_2, t_3, t_4) \\
&=&{}_4\phi_3\left( \left. \ba{c}
 q^{-n}, -t_1t_2q^{n-2}, -t_1t_3/q, -t_1t_4/q \\
 -t_1e^{\xi},\;  t_1e^{-\xi}, \; t_1t_2t_3t_4q^{-3} \ea \right| q,
q\right), 
\nonumber
\eea
\bea
\tilde{R}_n(z) = R_n(z; t_2, t_1, t_3, t_4), \nonumber
\eea 
and from (4.29) we have
\be
\frac{P_n(z)}{\prod_{j=1}^n(z - b_{j+1})}=
(-1)^nq^{n(n+1)/2}(t_3t_4)^{-n/2}
\frac{(-t_1t_2/q^2, t_1t_2t_3t_4q^{-3}; q)_{n}}
{(-t_1t_2/q^2; q)_{2n} \prod_{k=1}^n
u_k}\; R_n(z).
\ee
For the case $m=n$ we use (3.16) with $k = n$ and $z^n$ replaced by

$P_n(z; t_2, t_1, t_3, t_4)$ and $ z = a_{n+2}$. This requires the
evaluations
\be
P_n(a_{n+2}; t_2, t_1, t_3, t_4) = 
\frac{q^{-n(n-3)/2}(-t_2t_3/q, -t_2t_4/q;
q)_n}{(-2t_2\sqrt{t_3t_4})^n
\; \prod_{k=1}^n u_k}
\ee
and from  \cite [(III.24)]{Ga:Ra}
\be{}_8W_7(q^{n+1}t_4/t_3; q^{n+1}, -q^2/t_1t_3, t_4e^{-\xi},
-t_4e^{\xi},
-q^3/t_2t_3; t_1t_2/q)|_{e^\xi = t_2q^{n-1}}
\ee
\bea
= \frac{(q^{n+2}t_4/t_3, -t_1t_4/q, q^nt_2/t_3, -q^{2n-1}t_1t_2;
q)_\infty}
{(q^{2n+1}t_2/t_3,- q^nt_1 t_4, qt_4/t_3, -t_1t_2 q^{n-2};
q)_\infty}.
\nonumber
\eea
The final result after another lengthy calculation, is
\be
\sum_{k=0}^\infty R_m(z_k)\tilde{R}_n(z_k) \frac{(-q^2/t_1t_3,
-q^2/t_2t_3, -q^2/t_3t_4, -q^2/t_3^2; q)_k}{(qt_4/t_3, qt_1/t_3, q
t_2/t_3, q;q)_k}
\frac{1+q^{2k+2}/t_3^2}{1+q^2/t_3^2}(t_1t_2t_3t_4/q^3)^k
\ee
\bea
= (t_1t_2t_3t_4/q^3)^n\frac{1 + t_1t_2/q^2}{1 + t_1t_2q^{2n-2}}
\frac{(-q^2/t_3t_4, q; q)_n}{(-t_1t_2/q^2, t_1t_2t_3t_4/q^3; q)_n}
\nonumber
\eea
\bea
{}\qquad \times\frac{(-t_1t_2/q, -t_1t_4/q, -t_2t_4/q, -q^3/t_3^2;
q)_\infty}
{(qt_1/t_3, qt_2/t_3, qt_4/t_3, t_1t_2t_3t_4/q^3; q)_\infty}\;
\delta_{m,n},
\nonumber
\eea
with 
\bea
z_k = \frac{1}{2}(t_3q^{-k-1}- q^{k+1}/t_3), \; k = 0, 1, \cdots,
\quad |t_1t_2t_3t_4| < q^3.
\nonumber
\eea

Another biorthogonality relation can be obtained by 
interchanging $t_3$ and $t_4$, which corresponds to taking 
$u = -\sqrt{t_4/t_3}$. Both of these are special cases of the
biorthogonality
 relations derived in \cite {Is:Ma}, see (1.21), (3.14) and (3.15)
in \cite {Is:Ma}, with $a = q/t_3$ or $a = q/t_4$. However, if in
(4.28) we had used $\b = q 
e^{\xi}/\sqrt{t_3t_4}$ and $z = \cosh \xi$, then we would have
arrived 
at a similar biorthogonality with mass points at 
\bea
z_k =  \frac{1}{2}(t_3q^{-k-1} + q^{k+1}/t_3), k = 0, 1, \cdots.
\nonumber
\eea
We do not know where this fits into the scheme of things other than
to
 suggest that, with discrete orthogonalities at this level, there
are
 always two families, one associated with $z = \sinh \xi$, and
another with
$ z = \cosh \xi$.

\bigskip

\setcounter{section}{5}
 \setcounter{equation}{0}

{\bf 5. Biorthogonality on $[-1,1]$}.
We present two final examples by yet another modification of the
Askey-Wilson
 recurrence relation. Here an $R_{II}$-type  continued fraction is
used to 
obtain two $q$-beta integrals, one which is new and one which is in

\cite {Ra} . However the biorthogonalities 
obtained are not given by the continued fraction methods of the
 previous examples. Instead we use the ``attachment" method
explained in \cite {Be:Is} which goes back to Andrews and Askey and
was used in \cite {As:Wi}, 
\cite {Al:Is}. 

In (4.3) we apply the parameter replacements
\be
\b \to \b/u, \quad \gamma \to \b u,
\ee
with $z = (u+1/u)/2$. The recurrence relation again becomes an
$R_{II}$-type
recurrence relation (4.30) but now with
\be
\l_{n+1} = \frac{\a\b^2\delta q^{2n-2}(1-q^n)(1 - \a\b^2\delta
q^{n-2})
(1 - \a \delta q^{n-1})(1 - \b^2 q^{n-1})}
{(1 - \a \b^2\delta q^{2n-3})(1 - \a\b^2\delta q^{2n-2})^2
(1 - \a\b^2\delta q^{2n-1}) u_n u_{n+1}},
\ee
\be
a_{n+1} = \frac{1}{2}(\b\delta q^{n-1}+ q^{1-n}/\b \delta), \quad
b_{n+1} 
= \frac{1}{2}(\a\b q^{n-1} + q^{1-n}/\a\b), 
\ee
\be
c_{n+1} = v_{n+1}/u_{n+1},
\ee
where $u_n$ and $v_n$ are given by
\be
u_{n+1} = 1 - \b q^{n-1} \frac{(1 + \a\b^2\delta q^{2n-1})(q +
\a\delta)
- q^{n-1}(1+q)(q + \b^2)\a\delta }{(1 -\a\b^2\delta q^{2n-2})
(1 -\a\b^2\delta q^{2n})},
\ee
\be
v_{n+1} = (\a + \delta)q^{n-1} \frac{(1 + \a\b^2\delta q^{2n-1})(q
+ \b^2) -
q^{n-1}(1 + q)\b^2(q + \a\delta)}{2(1 - \a\b^2\delta q^{2n-2})
1 - \a\b^2\delta q^{2n})}.
\ee
The new minimal solution of the recurrence relation is given by
\bea
X_n^{(min)}(z) = \left\{ \ba{c}G_n(u), \quad |u| > 1
    \\ G_n(1/u), \quad |u| < 1, \ea \right.
 \eea
\bea
G_n(u) &=& \frac{(2u)^{-n}(\a\b^2\delta q^{2n-1}, \a q^{n+1}/u, \b 
q^{n+1}/u^2,
 \b q^{n+1}, \delta q^{n+1}/u; q)_\infty}
{(q^{n+1}, q^{n+2}/u^2, \a\b q^n u, \a\b q^n /u, \a\delta q^n,
\b^2q^n,
\b\delta q^n u, \b\delta q^n/u; q)_\infty} \\
&& \times {}_8W_7(q^{n+1}/u^2; q^{n+1}, q/\a u, q/ \b, qu^{-2}/ \b,
q/\delta u;
 \a\b^2\delta q^{n-1})/\prod_{k=-1}^n u_k. \nonumber
\eea

An application of Pincherle's theorem now yields the Stieltjes
transform 
\bea
\int_{-1}^1 \frac{f(x)\, dx}{z -x} &=& \frac{2}{u} \frac{(1 -
q/u^2)(1 - \a\b^2\delta /q)}{(1 - \a/u)(1 - \b/u^2)(1 - \b)(1
-\delta/u)}
\\
&{}& \mbox{}\qquad \times{}_8W_7(q/u^2; q, q/\a u, q/\b , q/\b u^2,
q/\delta u; \a\b^2\delta/q),
\nonumber
\eea
 with $z = \frac{1}{2}(u + 1/u), \; |u| > 1$,  and
\be
f(\cos \t) := \frac{1}{2\pi}  \frac{(e^{2i\t}, e^{-2i\t}, \a\b
e^{i\t}, 
\a\b  e^{-i\t}, \b\delta  e^{i\t}, \b\delta  e^{-i\t}, \a\delta,
 \b^2, q;q)_\infty}
{(\a  e^{i\t}, \a  e^{-i\t}, \b e^{2i\t} , \b e^{-2i\t}, \delta
e^{i\t}, \delta e^{-i\t}, \b, \b, \a\b^2\delta ; q)_\infty}
\frac{1}{\sin \t}.
\ee
The integral (5.9) seems to be new.

The large $z$ asymptotics of the
 transform formula (5.9) produces the following result
\bea
\frac{1}{2\pi} \int_{0}^\pi f(\cos \t) \, \sin \t  \, d\t
 = \frac{1 -\a\b^2\delta /q}{1 - \b} {}_2\phi_1(q, q/\b; 
q\b; q, \a\b^2\delta /q), 
\eea
valid for $\max\{|\a|, |\b|, |\delta|, |\a \b^2 \delta/q|
\} < 1$.

We next explain where (5.11) comes from in the theory of orthogonal
polynomials. 
Recall that the continuous  $q$-ultraspherical polynomials
$\{C_n(x;\b|q)\}$
 have the generating function
\bea
\sum_{n=0}^\infty C_n(\cos \t;\b|q) t^n =
 \frac{(\b te^{i\t}, \b te^{-i\t};q)_\infty}{( te^{i\t}, 
te^{-i\t};q)_\infty},
\eea
and satisfy the orthogonality relation
\bea
&&\frac{1}{2\pi} \int_0^\pi C_m(\cos \t;\b|q)C_n(\cos \t;\b|q) 
\frac{(e^{2i\t}, e^{-2i\t};q)_\infty}{(\b e^{2i\t}, \b
e^{-2i\t};q)_\infty} d\t \\
&& \mbox{} \qquad = \frac{(\b, q\b; q)_\infty(\b^2;q)_n(1-\b)}{(q,
\b^2;q)_\infty(q;q)_n(1-\b q^n)}
\; \delta_{m,n}, \nonumber
\eea
\cite{As:Is1}, \cite{Ga:Ra}. It is easy to see that (5.12) and
(5.13)  
show that 
the left-hand side of (5.11) is 
\bea
&=& \sum_{m,n=0}^\infty \frac{\a^n \delta^m}{2\pi} \int_0^\pi
C_m(\cos \t;\b|q)C_n(\cos \t;\b|q) 
\frac{(e^{2i\t}, e^{-2i\t};q)_\infty}{(\b e^{2i\t}, \b
e^{-2i\t};q)_\infty} d\t
\nonumber\\
&=& \sum_{n=0}^\infty (\a \delta)^n\frac{(\b, q\b;
q)_\infty(\b^2;q)_n(1-\b)}{(q, \b^2;q)_\infty(q;q)_n(1-\b q^n)}
\nonumber \\
&& \mbox{} \qquad = \frac{(\b, q\b;q)_\infty}{q, \b^2;q)_\infty}\;
{}_2\phi_1(\b^2, \b; q\b;q, \a\delta). \nonumber
\eea
Now the Heine transformation \cite[(III.3)]{Ga:Ra} reduces the
extreme
 right-hand side above to the right-hand side of (5.11).

\bigskip

{\bf Example 5.1 Chebyshev Rational Functions}. 
The relationship  (5.11) gives an integral representation for 
a  basic hypergeometric function
 of the type $_2\phi_1$. In this generality the integral
representation (5.11) 
does not seem to lead to
 orthogonal or biorthogonal functions. The special case $\b = q$
 degenerates into the elementary integral result 
(1.3).  To find rational functions  biorthogonal with respect to
the integrand in (1.3) we now use the attachment method. 

Consider  functions of the type
\be
g_n(\cos \t; \a, \delta) := \sum_{k=0}^n \frac{(q^{-n}, \a e^{i\t},
 \a e^{-i\t}; q)_k}{(q, q\a e^{i\t}, q\a e^{-i\t}; q)_k}a_{n,k},
\ee
where $a_{n.k}$ are to be determined.  Set
\be
I_{n,j} = \frac{2}{\pi}\int_{0}^\pi  g_n(\cos \t; \a, \delta) 
 \frac{(\delta e^{i\t}, 
\delta e^{-i\t}; q)_j}{(q\delta e^{i\t}, q\delta e^{-i\t}; q)_j} 
\frac{\sin^2 \t d\t}{(1- 2\a \cos \t + \a^2)(1 - 2\delta \cos \t +
\delta^2)}.
\ee
Thus 
\bea
I_{n,j} &=& \sum_{k=0}^n \frac{(q^{-n}; q)_k}{(q; q)_k} a_{n,k} 
\frac{2}{\pi} \int_{-1}^1 \frac{\sqrt{1-x^2}\, dx}{(1 - 2 \a q^k x
+ \a^2 q^{2k})(1 - 2\delta q^j x +\delta^2q^{2j})} \nonumber\\
 &=&   \sum_{k=0}^n \frac{(q^{-n}; q)_k}{(q; q)_k} a_{n,k}
\frac{1}{1 - \a\delta q^{k+j}} \nonumber \\
&=& \frac{1}{1 - \a\delta q^j}  \sum_{k=0}^n 
\frac{(q^{-n}, \a\delta q^j; q)_k}{(q, \a\delta q^{j+1}; q)_k}
a_{n,k}.
\nonumber
\eea
Recall that  the $q$-analogue of the Pfaff-Saalsch\"utz theorem
\cite [II.12)]
{Ga:Ra} is
\be
{}_3\phi_2\left( \left. \ba{c}
 q^{-n}, a, b\\
 c, abq^{1-n}/c \ea \right| q, q\right) = \frac{(c/a, c/b; q)_n}
{(c, c/ab; q)_n}.
\ee
We would like to choose $a_{n,k}$ to make $I_{n,j}$ vanish for $0
\le j < n$. 
In view of (5.16) we let $a_{n,k} = q^k(\ a\delta
q^n;q)_k/(\a\delta; q)_k$ and
we obtain
\bea
I_{n,j} = \frac{(q^{-n}, q^{-j}; q)_n}{(1-\a\delta q^j)\,
(\a\delta, q^{-n-j}/\a\delta)_n}.  \nonumber
\eea
It is now clear that $I_{n,j} = 0$ if $0 \le j < n$. After some
simplification we find
\bea
I_{n,j} = [(-\a\delta)^nq^{n(n-1)/2}(q; q)^2_n/
(\a\delta; q)_{2n+1}] \delta_{j,n}. \nonumber
\eea
Thus  the $g_n$'s of (5.14) are given by
\be
g_n(\cos \t; \a, \delta) = {}_4\phi_3\left( \left. \ba{c}
 q^{-n}, \a\delta q^{n}, \a e^{i\t}, \a e^{-i\t}\\
 \a\delta,\;  q\a e^{i\t}, \; q\a e^{-i\t} \ea \right| q, q\right)
\ee
and satisfy the orthogonality relation
\be
\frac{2}{\pi} \int_{-1}^1 g_m(x; \delta,  \a) g_n(x; \a, \delta)
\frac{\sqrt{1-x^2}\, dx}{(1 - 2 \a x + \a^2)(1 - 2\delta x +
\delta^2)}
\ee
\bea
{}\qquad \qquad  = \frac{(\a\delta)^n\, (q, q; q)_n}{(\a\delta,
\a\delta ;q)_n
(1 - \a\delta q^{2n})} \;\delta_{m,n}. \nonumber
\eea

We find the biorthogonality relation (5.18) very surprising since
the 
${}_4\phi_3$ function in (5.17) is not balanced. Recall that  a
basic 
hypergeometric function (1.9) is balanced if $r = s+1$ 
 and $qa_1 a_2  \cdots a_{s+1} = b_1b_2 \cdots b_s$. Only balanced 
${}_4\phi_3$'s  with argument $q$ satisfy three term recurrence
relations.  
 This concludes this example.

\bigskip

 In the previous examples of pairs
of biorthogonal functions the second family was obtained from the 
 first by symmetry and permutation of  parameters.  In the next
example the 
members of the pair of biorthogonal rational functions are not
related
in this manner.

\bigskip

{\bf Example 5.2  A Nonsymmetric Case}.  Following  previous
examples
 we now consider the $q$-beta integral obtained by absorbing the
factor
 $x - a_2$ into
the weight function. This means evaluating (5.9) at $z = a_2$, that
is 
$u = 1/\a \b$, and this is exactly the case that makes the
right-hand
side of (5.9) summable. This gives the  result
\bea
 \frac{1}{2\pi} \int_{0}^\pi \frac{(e^{2i\t}, e^{-2i\t}, \a\b q
e^{i\t}, 
\a\b  qe^{-i\t},
\b\delta e^{i\t}, \b\delta e^{-i\t}, \a\delta, \b^2, q; q)_\infty}
{(\a e^{i\t}, \a e^{-i\t},\b e^{2i\t}, \b e^{-2i\t},\delta e^{i\t},
\delta
 e^{-i\t}, \b , q\b , \a\b^2\delta; q)_\infty} \; d\t\\
{}\qquad = \frac{1}{1 - \a^2 \b}, \qquad max\{|\a|, |\b|,
|\delta|\} < 1.
\nonumber
\eea

After we showed an earlier version of this paper to Mizan Rahman he
pointed out that  this $q$-beta integral is not new. It was
obtained
previously by him in his interesting work \cite {Ra} using a
different method. 
It may be of interest to mention here that  (5.19) is also the
special case
$\g = q\b$  of
\bea
&& \frac{1}{2\pi} \int_{0}^\pi \frac{(e^{2i\t}, e^{-2i\t}, \a\g
e^{i\t}, 
\a\g e^{-i\t},
\b\delta e^{i\t}, \b\delta e^{-i\t}, \a\delta, \b^2, q; q)_\infty}
{(\a e^{i\t}, \a e^{-i\t},\b e^{2i\t}, \b e^{-2i\t},\delta e^{i\t},
\delta
 e^{-i\t}, \b , q\b , \a\b^2\delta; q)_\infty} \; d\t \\
&&{}\qquad = \frac{(\b, \g,\a^2\g, \b^2\a\delta;q)_\infty}{(q, 
\b^2, \a^2\b,\a\delta;q)_\infty}{}_3\phi_2\left( \left. \ba{c}
  \a^2\b, \a\delta, q\b/\g\\
 \a^2\g, \b^2\a\delta \ea \right| q, \g\right), 
\nonumber
\eea
 for  max$\{|\a|, |\b|, |\g|,
 |\delta|\} < 1$, which will appear elsewhere. In view of (5.12)
and (5.13) formula (5.20) 
is equivalent to the evaluation 
of the connection coefficients for the continuous
$q$-ultraspherical polynomials.
 L.J. Rogers solved this connection coefficient problem in 1893,
see \cite{As:Is1}.

One can use
the attachment method of the previous example to derive a 
system  biorthogonal with respect to the integrand in (5.19). 
Rahman also has done this in \cite{Ra}. For completeness we repeat
the
calculation here.

The $q$-beta integral we are working with is
\be
\int_{-1}^1  w(x; \a, \b, \delta) dx
 = \frac{(\b, q\b, \a\b^2\delta; q)_\infty}{(1 - \a^2\b)(\a\delta,
\b^2, q; q)_\infty},
\ee
where
\be
w(\cos \t; \a, \b, \delta) := \frac{(e^{2i\t}, e^{-2i\t}, q\a\b
e^{i\t}, 
q\a\b e^{-i\t},
\b\delta e^{i\t}, \b\delta e^{-i\t}; q)_\infty}{2\pi\sin \t \,(\a
e^{i\t}, \a e^{-i\t},\b e^{2i\t}, \b e^{-2i\t},\delta e^{i\t},
\delta
 e^{-i\t}; q)_\infty}.
\ee
Choose 
\be
\psi_n(\cos \t; \a, \b, \delta) = \sum_{k=0}^n \frac{(q^{-n},
\delta e^{i\t},
 \delta e^{-i\t}; q)_k}{(q, \delta\b e^{i\t}, \delta \b e^{-i\t};
q)_k}a_{n,k}.
\ee
Using the   $q$-analogue of the Pfaff-Saalsch\"utz   theorem,
(5.16) we put
\bea
a_{n,k} = q^k\frac{(\a\b^2\delta q^{n-1}; q)_k}{(\a\delta; q)_k},
\nonumber
\eea
and find 
\bea
 \int_{0}^\pi w(\cos \t; \a, \b, \delta) 
\psi_n(\cos \t; \a, \b, \delta)  \frac{(\a e^{i\t}, 
\a e^{-i\t}; q)_j}{(q\a\b e^{i\t}, q\a\b e^{-i\t}; q)_j} \sin \t
d\t. \\
=  \frac{(\b, q\b, \a\b^2q^{n}\delta; q)_\infty}{(1 - \a^2\b
q^{2n})
(\a\delta q^n, \b^2, q;q)_\infty} \frac{(\b^2, q;
q)_n}{(\a\b^2\delta q^{n}, 
q^{1-n}/\a\delta; q)_n} \delta_{j, n}, \quad j \le n.
 \nonumber
\eea
Let
 \be
\psi_n(\cos \t; \a, \b, \delta) := {}_4\phi_3\left( \left. \ba{c}
 q^{-n}, \a\b^2\delta q^{n-1}, \delta e^{i\t}, \delta e^{-i\t}\\
 \a\delta,\;  \b\delta e^{i\t}, \;  \b\delta e^{-i\t} \ea \right|
q, q\right), 
\ee
and
\be
 \phi_n(\cos \t; \a, \b, \delta):= \sum_{j=0}^n 
\frac{(q^{-n}, \a e^{i\t}, \a e^{-i \t}; q)_j}
{(a, a\a\b e^{i\t}, q\a\b e^{-i\t}; q)_j} q^j b_{n,j}.
\ee
Thus we have proved that the rational functions $\{\phi_n\}$ and 
$\{\psi_n\}$ 
satisfy the biorthogonality relation
\bea
\int_{-1}^1 \phi_m(x; \a , \b , \delta) \psi_n(x; \a, \b , \delta)
w(x; \a , \b , \delta) dx  
\eea
\bea
 \qquad = \frac{(\b, q\b, \a\b^2\delta q^{n}; q)_\infty}{(\a\delta
, 
\b^2, q;q)_\infty} \frac{(\b^2, q; q)_n}
{(1 - \a^2\b q^{n})} (-\a\delta)^n \; b_{n,n}\;\delta_{m,n}, \quad
m \le n.\nonumber
\eea
We now select $b_{n,j}$ in order to have full biorthogonality, that
is 
to make (5.27) hold for all $m$ and $n$.
 In order to do so we need to compute the integrals $J_{n,k}$,
\bea
J_{n,k} := \int_{0}^\pi \phi_n(\cos \t; \a , \b , \delta) 
\frac{(\delta e^{i\t}, \delta e^{-i\t}; q)_k}
{(\b\delta e^{i\t},\b \delta e^{-i\t}; q)_k}
w(\cos \t; \a , \b , \delta) \sin \t d\t, \nonumber
\eea
for $k\le n$. Substituting for $\phi_n$ from (5.26) and using 
the $q$-beta integral (5.21) we get
\bea
J_{n,k} =  \frac{(\b, q\b, \a\b^2\delta q^{k}; q)_\infty}{(\b^2,
\a\delta q^k
; q)_\infty} \sum_{j = 0}^n \frac{(q^{-n}, \a\delta q^k; q)_j}{q,
\a\b^2 \delta 
q^k; q)_j}\; \frac{q^j\,b_{n,j}}{1 - \a^2\b q^{2j}}.
\nonumber
\eea
In order to apply (5.16) we must choose $b_{n,j}$ as
\be
b_{n,j} = \frac{(1 - \a^2\b q^{2j})(\a\b^2\delta q^{n-1}; q)_j}{(1
- \a^2\b)\;(\a\delta;
q)_j},
\ee
hence 
\bea
J_{n,k} = \frac{(\b, q\b, \a\b^2\delta q^{k}; q)_\infty}{(\b^2,
\a\delta q^k
; q)_\infty} \frac{(q^{-k}, q^{1-n}/\b^2; q)_n}
{(\a\delta, q^{1-k-n}/a\b^2\delta; q)_\infty}. \nonumber 
\eea
Thus (5.26) becomes
\be
\phi_n(\cos \t, \a, \b, \delta) = {}_6\phi_5\left( \left. \ba{c}
 q^{-n}, \a\b^2\delta q^{n-1}, \a e^{i\t}, \a e^{-i\t},
q\a\sqrt{\b}, 
- q\a\sqrt{\b}\\
 \a\delta,\; q \a\b e^{i\t}, \;  q\a\b e^{-i\t}, \a\sqrt{\b}, 
- \a\sqrt{\b}\ea \right| q, q\right).
\ee

This establishes the biorthogonality relation
\bea
\int_{-1}^1 \phi_m(x; \a , \b , \delta) \psi_n(x; \a, \b , \delta)
w(x; \a , \b , \delta) dx 
\eea
\bea
\qquad = 
 \frac{(\b, q\b, \a\b^2 \delta q^{n}; q)_\infty}
{(1-\a^2\b)(\a\delta , \b^2, q;q)_\infty}
 \frac{(\b^2 , q, ; q)_n\; (1 - \a^2\b \delta q^{2n-1}) \;
(\a\delta)^n}
{(1 - \a^2\b \delta q^{n-1}) 
(\a\delta; q)_n}  \delta_{m,n}, \nonumber
\eea
valid for $|\a|, |\b|, |\delta| < 1$.

Special cases of the $\psi_n$'s have appeared in the literature as 
special cases of the Askey-Wilson polynomials. The case $\b = 0$
 gives the Al-Salam-Chihara 
polynomials, \cite {As:Is2}, \cite {As:Wi}, as can be seen by 
comparing (5.26)
 and (5.29) with (4.1). The case $\a = \delta = 0$ of  (5.21),
 (5.22), (5.26) and (5.29) gives the $q$-ultraspherical 
polynomials of L. J. Rogers \cite {As:Is1}, \cite {As:Wi}
 but the justification is not as straightforward. It can be
justified however by using the generating function \cite {Is:Wi}
\bea
\sum_0^\infty \frac{(ac, ad;q)_n}{(q,cd; q)_n}(t/a)^n 
{}_4\phi_3\left( \left. \ba{c}
 q^{-n}, abcd q^{n-1}, a e^{i\t}, a e^{-i\t}   \\
 ab, \; ac, \; ad \ea \right| q, q\right)
\eea
\bea
= {}_2\phi_1\left( \left. \ba{c}
  a e^{i\t}, b e^{i\t}   \\
 ab \ea \right| q, q\right)\; {}_2\phi_1\left( \left. \ba{c}
  c e^{-i\t}, d e^{-i\t}   \\
 cd \ea \right| q, q\right).   \nonumber
\eea
In a private communication R. Askey mentioned that he can let $\a
\to 0$ 
in (4.1) and show directly that 
$\lim_{\a \to 0} \a^{-n}P_n(x; \a, \b, \gamma, \delta|q)$ exists
and
 find its value. We do 
not know how to compute the limit 
$\lim_{\a \to 0} \a^{-n}\phi_n(x; \a, \b, \delta|q)$ directly  but
 we suspect that Askey's direct argument may work here. In \cite
{Ra}
 M. Rahman made the transition to the $\a=\delta=0$ case by 
 using a ${}_3\phi_2$ representation
 for the $q$-ultraspherical polynomials, writing $\phi_n$ as a
 combination of ${}_4\phi_3$'s, applying Sear's transformation for
 ${}_4\phi_3$'s to the $\psi_n$ and the ${}_4\phi_3$'s in $\phi_n$,
 and then taking the limit $\alpha, \delta \to 0 $.

\bigskip
\noindent{\bf Acknowledgments}. We thank the referee for his/her
careful 
and critical  reading of
the first version of this paper and for many helpful suggestions.
We also thank 
our friend Mizan Rahman for comments and for pointing out the
overlap of 
Example 5.2 and his  earlier work \cite{Ra}.

Department of Mathematics, University of South Florida, Tampa, 
Florida,
  USA 33620.

 Department of Mathematics, University of Toronto, Toronto,
  Ontario, Canada M5S 1A1

\bigskip

\begin{thebibliography}{99}
\bibitem{Al:Is} W. A. Al-Salam and M. E. H. Ismail,
 {\em A $q$-beta integral on the unit circle and some biorthogonal
 rational functions},
Proc. Amer. Math. Soc. {\bf 121}(1994), 553-561.
\bibitem{As}R. A. Askey, Comments on [21-1], in "Gabor Szego:
Collected Papers", editor R. A. Askey, volume 1, Birkhauser,
Basel, 1982, pp. 303-305.
\bibitem{As:Is1} R.~A.~Askey and M.~E.~H.~Ismail, {\it A
generalization
of ultraspherical polynomials}, in ``Studies in Pure Mathematics",
ed.~P.~Erd\"{o}s,
Birkhauser, Basel, 1983, pp.~55-78.
\bibitem{As:Is2} R.~A.~Askey and M.~E.~H.~Ismail, {\it Recurrence
relations, continued fractions and orthogonal polynomials}, Memoirs
Amer.~Math.~Soc. Number {\bf 300} (1984).
\bibitem{As:Wi}R. A. Askey and J. A. Wilson, {\em Some basic
hypergeometric polynomials that generalize Jacobi polynomials},
Memoirs Amer. Math. Soc. Number {\bf 319}, 1985.
\bibitem{Be:Is}C. Berg and M. E. H. Ismail, {\em $Q$-Hermite
polynomials and classical orthogonal polynomials}, to appear.
\bibitem{Ch}T. S. Chihara, An Introduction to Orthogonal 
Polynomials, Gordon and Breach, New York, 1978.
\bibitem{Er:Ma1}A. Erdelyi, W. Magnus, F. Oberhettinger 
and F. G. Tricomi,
Higher Transcendental Functions, Volume 1, McGraw-Hill, New York,
1953.
\bibitem{Ga:Ra}G. Gasper and M. Rahman, 
Basic Hypergeometric Series, Cambridge University Press, 
Cambridge, 1990. 
\bibitem{Go} A. A. Goncar, {\em On the speed of rational
approximation
of some analytic functions}, Math. USSR Sbornik {\bf 34}(1978),
131-145.
\bibitem{Go:Lo} A. A. Goncar and G. Lopes, {\em On Markov's theorem
for
multipoint Pad\'e approximants}, Math. USSR Sbornik {\bf 34}
(1978) 449-459.
\bibitem{Gu:Is}D. P. Gupta, M. E. H. Ismail and D. R. Masson, {\em 
Contiguous relations, basic hypergeometric series and orthogonal
polynomials
II: Associated big $q$-Jacobi polynomials},
{\mg J. Math. Anal. Appl.} {\bf 171} 
(1992), 477-497.
\bibitem{Gu:Ma}D. P. Gupta and D. R. Masson, {\em Solutions to the
q-Askey-Wilson polynomial recurrence relation}, to appear.
\bibitem{He:Nj} E. Hendriksen and O. Njasted, {\em A Favard theorem
for rational functions}, J. Math. Anal. Appl. {\bf 142} (1989),
508-520.
\bibitem{He:Va}E. Hendriksen and H. van Rossum, {\em Orthogonal
Laurent 
polynomials}, Indag. Math. (ser. A) {\bf 89} (1986), 17-36.
\bibitem{Is:Le:Va:Wi}M. E. H. Ismail, J. Letessier, G. Valent and
J. Wimp, {\em Two families of associated Wilson polynomials},
Can. J. Math. {\bf 42} (1990), 659-695.
\bibitem{Is:Li}M. E. H. Ismail and C. A. Libis, {\em Contiguous 
relations
basic hypergeometric functions and orthogonal polynomials I},
J. Math. Anal.  Appl. {\bf 141} (1989), 349-372.
\bibitem{Is:Ma}M. E. H. Ismail and D.  R. Masson, {\em Q-Hermite
polynomials, biorthogonal rational functions and $q$-beta
integrals},
Trans. Amer. Math.Soc. (1995), to appear.
\bibitem{Is:Wi}M. E. H. Ismail and J. A. Wilson, {\em Asymptotic
and generating relations for the $q$-Jacobi polynomials and the
${}_4\phi_3$ polynomials},
J. Approx. Theory {\bf36} (1982), 43-54.
\bibitem{Jo:Nj:Th}W. B. Jones, O. Njastad  and W. J. Thron,
{\em Moment theory, orthogonal polynomials, and continued fractions
associated with the unit circle}, Bull. London Math. Soc.
{\bf 21}  (1989),
113-152.
\bibitem{Jo:Th1}W. B. Jones and W. Thron, Continued Fractions:
Analytic Theory and Applications, Cambridge University Press,
Cambridge, 1980.
\bibitem{Lo:Wa}L. Lorentzen and H. Waadeland, Continued Fractions,
North-Holland, Amsterdam, 1992.
\bibitem{Ma}D. R. Masson, {\em Associated Wilson polynomials}, 
Constructive Approximation
{\bf 7} (1991), 521-534.
\bibitem{Ma2}D. R. Masson, {\em Explicit spectral theory: from
Chebyshev to Askey-Wilson}, to appear.
\bibitem{Nj}O. Njastad, {\em  Multipoint Pad\'e
approximation and orthogonal
rational functions} in "Nonlinear Numerical Methods and Rational 
Approximation", editor A. Cuyt, D. Reidel, Dordrecht 1988, pp. 
259-270.
\bibitem{Pa}P. I. Pastro, {\em Orthogonal polynomials and some
$q$-beta integrals of Ramanujan}, J. Math. Anal. Appl. {\bf 112} 
(1982), 517-540.
\bibitem{Ra}M. Rahman, {\it Some extensions of
Askey-Wilsons q-beta integral and the corresponding orthogonal
systems}, Canad. Math. Bull. {\bf 31} (1988), 111-120.
\bibitem{St:To}H. Stahl and V. Totik, General Orthogonal
Polynomials,
 Cambridge University Press, Cambridge, 1992.
\bibitem{Sz}G. Szeg\"o, Othogonal Polynomials, fourth edition,
American
Mathematical Society, Providence, 1975.
\bibitem{Th:Wa}W. Thron and H. Waadeland,
{\em On a certain transformation of continued  fractions} in
"Analytic Theory of Continued Fractions", eds. W. B. Jones, W.   
Thron and  H. Waadeland, Lecture Notes in Mathematics, number
932, Springer-Verlas, Berlin, 1981, pp. 225-240.
\bibitem{Wa}H. S. Wall, Analytic Theory of Continued Fractions, D.
Van Nostrand, Princeton, 1948.
\bibitem{Wi1}J. Wimp, {\em Some explicit Pad\'e approximants
for the function
$\phi'/\phi$ and a related quadrature formula involving 
Bessel functions}, SIAM J. Math. Anal. {\bf 16} (1985),
pp. 887-895.
\bibitem{Wi2}J. Wimp, {\em Explicit formulas for the associated
Jacobi polynomials and some applications}, Can. J. Math. {\bf 39}
 (1987), 983-1000.
\end{thebibliography}
\end{document}